\documentclass{iopart}
\usepackage{iopams,graphicx,latexsym,pdfsync}

\newcommand{\bigO}{\mathcal O}

\newcommand{\dk}{\, \mathrm{d}k}
\newcommand{\ds}{\, \mathrm{d}s}
\newcommand{\dx}{\, \mathrm{d}x}

\newcommand{\nn}{|{\mskip-2mu}|{\mskip-2mu}|}

\newcommand{\Efunc}{\mathcal{E}}
\newcommand{\Gfunc}{\mathcal{G}}
\newcommand{\Nfunc}{\mathcal{N}}
\newcommand{\Lfunc}{\mathcal{L}}
\newcommand{\Qfunc}{\mathcal{Q}}

\newcommand{\R}{\mathbb{R}}
\newcommand{\Z}{\mathbb{Z}}
\newcommand{\N}{\mathbb{N}}

\newcommand{\FF}{\mathcal{F}}
\renewcommand{\SS}{\mathcal{S}}

\newcommand{\dist}{\mathop{\mathrm{dist}}}

\newcommand{\supp}{\mathop{\mathrm{supp}}}

\newcommand{\proof}{{\bf Proof }}

\newcommand{\qed}{\hfill$\Box$\smallskip}

\newcommand{\tfrac}[2]{{\textstyle\frac{#1}{#2}}}

\newcommand{\afl}{\fl\quad}

\newtheorem{theorem}{Theorem}[section]
\newtheorem{lemma}[theorem]{Lemma}
\newtheorem{corollary}[theorem]{Corollary}
\newtheorem{proposition}[theorem]{Proposition}

\newtheorem{remark}[theorem]{Remark}

\begin{document}

\title[Solitary waves for equations of Whitham type]{On the existence and stability of solitary-wave solutions to a class of evolution equations of Whitham type}
\author{Mats Ehrnstr\"om}
\address{Department of Mathematical Sciences, Norwegian University of Science and Technology, 7491 Trondheim, Norway}
\address{Institut f\"ur Angewandte Mathematik, Leibniz Universit\"at Hannover, Welfengarten 1, 30167 Hannover, Germany}
\ead{mats.ehrnstrom@math.ntnu.no}
\author{Mark D. Groves}
\address{Fachrichtung 6.1 - Mathematik, Universit\"at des Saarlandes, Postfach 151150, 66041 Saarbr\"ucken, Germany}
\address{Department of Mathematical Sciences, Loughborough University, Loughborough, LE11 3TU, UK}
\ead{groves@math.uni-sb.de}
\author{Erik Wahl\'en}
\address{Centre for Mathematical Sciences, Lund University, PO Box 118, 221\,00 Lund, Sweden}
\ead{erik.wahlen@math.lu.se}

\begin{abstract}
We consider a class of pseudodifferential evolution equations of the form
$$u_t + ( n(u) + Lu )_x = 0,$$
in which $L$ is a linear smoothing operator and $n$ is at least quadratic near the origin;
this  class includes in particular the Whitham equation. A family of solitary-wave solutions
is found using a constrained minimisation principle and concentration-compactness methods
for noncoercive functionals. The solitary waves are approximated by
(scalings of) the corresponding solutions to partial differential equations arising
as weakly nonlinear approximations; in the case of the Whitham equation the approximation
is the Korteweg-deVries equation. We also demonstrate that the family of
solitary-wave solutions is conditionally energetically stable.
\end{abstract}
\ams{35Q53, 35A15, 76B15}
\maketitle




\section{Introduction}\label{sec:intro}

In this paper we discuss solitary-wave solutions of the pseudodifferential equation
\begin{equation}\label{eq:DEproblem}
u_t + (Lu+n(u))_x = 0
\end{equation}
describing the evolution of a real-valued function $u$ of time $t \in \R$ and space $x \in \R$;
here $L$ is a linear smoothing operator and $n$ is at least quadratic near the origin.
A concrete example is the equation
\begin{equation}\label{eq:whitham}
u_t + 2 u u_x + Lu_x = 0,
\end{equation}
where $L$ is the spatial Fourier multiplier operator given by
$$
\FF(Lf)(k) = \left( \frac{\tanh(k)}{k} \right)^{\!\!\frac{1}{2}} \hat f(k).
$$
This equation was proposed by Whitham~\cite{Whitham67} as an alternative to the
Korteweg-deVries equation which features the same linear dispersion relation as the full water-wave problem,
a fact that allows for the breaking of waves (Whitham~\cite{Whitham}, Naumkin \& Shishmarev~\cite{NaumkinShishmarev}).
There have been
several investigations of different variants of the Whitham equation (e.g.\ see Constantin \& Escher~\cite{ConstantinEscher98},
Gabov~\cite{Gabov78}, Naumkin \& Shishmarev~\cite{NaumkinShishmarev} and Zaitsev~\cite{Zaitsev86}), but it has remained
unclear whether the Whitham equation admits travelling waves, that is solutions of the form $u=u(x-\nu t)$ representing
waves moving from left to right with constant speed $\nu$.
The existence of periodic travelling waves to the Whitham equation
was recently established by Ehrnstr\"{o}m \& Kalisch~\cite{EhrnstroemKalisch09},
and in the present paper we discuss solitary waves, that is travelling waves for which $u(x-\nu t) \rightarrow 0$
as $x-\nu t \rightarrow \pm\infty$.

Our mathematical task is therefore to find functions $u=u(x)$ which satisfy the travelling-wave equation
\begin{equation}\label{eq:steadyproblem}
Lu -\nu u + n(u) = 0
\end{equation}
with wave speed $\nu$ and asymptotic condition $u(x) \rightarrow 0$ as $x \rightarrow \pm \infty$. 
We examine equation~\eref{eq:steadyproblem} under the following conditions.\\
\\
{\bf Assumptions}{\it
\begin{itemize}
\item[(A1)]
The operator $L$ is a Fourier multiplier with classical symbol $m \in S_\infty^{m_0}(\R)$ for some $m_0 < 0$,
that is
$$
\FF(Lf)(k) = m(k) \hat f(k)
$$
for some smooth function $m \colon \R \to \R$ with the property that
\begin{equation}\label{eq:mclassical}
|m^{(\alpha)}(k)| \leq C_\alpha \left( 1 + |k| \right)^{m_0 - \alpha}, \qquad \alpha \in \N_0,
\end{equation}
where $C_\alpha$ is a positive constant depending upon $\alpha$. In particular, one can write
$L$ as a convolution with the (possibly distributional) kernel $K := \FF^{-1}(m)$, that is
\begin{equation}\label{eq:Lconvolution}
Lf = \frac{1}{\sqrt{2\pi}}\, K \ast f.
\end{equation}
\item[(A2)]
The symbol $m: \R \to \R$ is even (to avoid non-real solutions) and satisfies $m(0)>0$,
\begin{equation}\label{eq:mmax}
m(k) < m(0), \qquad k \neq 0,
\end{equation}
(so that it has a strict and positive global maximum at $k = 0$) and
$$
m(k) = m(0) + \frac{m^{(2j_\star)}(0)}{(2j_\star)!}\, k^{2j_\star} + r(k)
$$
for some $j_\star \in \N$, where $m^{(2j_\star)}(0) < 0$ and $r(k)=\bigO(k^{2j_\star+2})$ as $k \to 0$.
\item[(A3)]
The nonlinearity $n$ is a twice continuously differentiable function $\R \rightarrow \R$
with
\begin{equation}\label{eq:nbigO}
n(x) = n_p(x) + n_\mathrm{r}(x),
\end{equation}
in which the leading-order part of the nonlinearity takes the form
$n_p(x)=c_p |x|^p$ for some $c_p \neq 0$ and $p \in [2,4j_\star+1)$
or $n_p(x)=c_p x^p$ for some $c_p>0$ and odd integer $p$ 
in the range $p \in [2,4j_\star+1)$, while the higher-order part of the nonlinearity satisfies the estimate
$$n_\mathrm{r}(x) = \bigO(|x|^{p+\delta}), \qquad n_\mathrm{r}^\prime(x)=\bigO(|x|^{p+\delta-1})$$
for some $\delta>0$ as $x \to 0$. (Occasionally we simply estimate $n(x) = \bigO(|x|^p)$ and
$n^\prime(x)=\bigO(|x|^{p-1})$.)
\end{itemize}}

Proceeding formally, let us derive a long-wave approximation to equation~\eref{eq:steadyproblem}
by introducing a small parameter $\mu$ equal to the momentum $\frac{1}{2}\int_\R u^2\dx$ of the wave,
writing $\nu$ as a small perturbation of the speed $m(0)$
of linear long waves, so that
$$\nu=m(0) + \mu^\gamma \nu_\mathrm{lw},$$
and substituting the weakly nonlinear \emph{Ansatz}
\begin{equation} \label{eq:longwaveansatz}
u(x) := \mu^{\alpha} w(\mu^{\beta} x),
\end{equation}
where $2\alpha-\beta=1$ (so that $\frac{1}{2}\int_\R u^2=\mu$) into the equation.
Choosing $(p-1)\alpha=2j_\star\beta$ and $\gamma=2j_\star\beta$, we find that
$$\mu^{p\alpha}\left(\frac{(-1)^{j_\star}}{(2j_\star)!}
m^{(2j_\star)}(0)w^{(2j_\star)} - \nu_\mathrm{lw} w + n_p(w)\right) + \ldots=0,$$
where the ellipsis denotes terms which are formally $o(\mu^{p\alpha})$; the constraints on $\alpha$
and $\beta$ imply the choice
\begin{equation}\label{eq:alphabeta}
\alpha = \frac{2j_\star}{4j_\star +1-p} \quad\mbox{ and }\quad  \beta = \frac{p-1}{4j_\star +1-p}.
\end{equation}

This formal
weakly nonlinear analysis suggests that solitary-wave solutions to~\eref{eq:DEproblem} are
approximated by (suitably scaled) homoclinic solutions of the ordinary differential equation
\begin{equation}\label{eq:lwe}
\frac{(-1)^{j_\star}}{(2j_\star)!}m^{(2j_\star)}(0)w^{(2j_\star)}-\nu_\mathrm{lw} w + n_p(w)=0
\end{equation}
for some constant $\nu_\mathrm{lw}$. The following theorem gives a variational characterisation
of such solutions; it is established using a straightforward modification of the theory developed
by Albert \cite{Albert99} and Zeng \cite{Zeng03} for a slightly different class of equations
(the proof that $\Efunc_\mathrm{lw}$ is bounded below over $W_1$ is given in the appendix).\pagebreak

\begin{theorem}\label{theorem:AlbertZeng}\hspace{1in}
\begin{itemize}
\item[(i)]
The functional
$\Efunc_\mathrm{lw}: H^{j_\star}(\R) \rightarrow \R$ given by
\begin{equation} \label{eq:redEdefn}
\Efunc_\mathrm{lw}(w) = -\int_\R\left\{\frac{m^{(2j_\star)}(0)}{2(2j_\star)!} (w^{(j_\star)})^2 +N_{p+1}(w)\right\}\dx,
\end{equation}
where
$$
N_{p+1}(x) := 
\left\{\begin{array}{ll}
\displaystyle \frac{c_p x^{p+1}}{p+1}, & \quad\mbox{if $n_p(x)=c_p x^p$},\\[5mm]
\displaystyle \frac{c_p x|x|^p}{p+1}, & \quad\mbox{if $n_p(x)=c_p |x|^p$},
\end{array}\right.
$$
is bounded below over the set
$$W_1=\{w \in H^{j_\star}(\R): \Qfunc(w)=1\},$$
where
\begin{equation}\label{eq:Qdefn}
\Qfunc(w) = \frac{1}{2}\int_\R w^2\dx.
\end{equation}
The set $D_\mathrm{lw}$ of minimisers of $\Efunc_\mathrm{lw}$ over $W_1$
is a nonempty subset of $H^{2j_\star}(\R)$ which lies in
$$W:=\{w \in H^{2j_\star}(\R): \|w\|_{2j_\star} < S\}$$
for some $S>0$.
Each element of $D_\mathrm{lw}$ is a solution of equation~\eref{eq:lwe};
the constant $\nu_\mathrm{lw}$ is the Lagrange multiplier in this constrained variational principle.
\item[(ii)]
Suppose that $\{w_n\}_{n \in \N_0}$ is a minimising sequence for $\Efunc_\mathrm{lw}$ over 
$\{w \in H^{j_\star}(\R): \Qfunc(w)=1\}$. There exists a sequence $\{x_n\}_{n \in \N_0}$ of real numbers with the property that a subsequence of $\{w_n(\cdot + x_n)\}_{n \in \N_0}$ converges in $H^{j_\star}(\R)$ to an element of
$D_\mathrm{lw}$.
\end{itemize}
\end{theorem}

For the Whitham equation ($j_\star=1$, $p=2$, $m^{\prime\prime}(0)=-\frac{1}{3}$)
the above derivation yields the travelling-wave version
$$\tfrac{1}{6}w^{\prime\prime} - \nu_\mathrm{lw} w + w^2=0$$
of the Korteweg-deVries equation,
for which
$$D_\mathrm{lw}= \{w_\mathrm{KdV}(\cdot + y): y \in \R\}, \qquad
w_\mathrm{KdV}(x) = \left(\tfrac{3}{2}\right)^{\!\frac{2}{3}}\mathrm{sech}^2 \left( \left(\tfrac{3}{2}\right)^{\!\frac{1}{3}} x \right)$$
(and there are no further homoclinic solutions). In general $D_\mathrm{lw}$ consists of all spatial translations of a (possibly
infinite) family of `generating' homoclinic solutions with different wave speeds
($\nu_\mathrm{lw} = \left(\frac{2}{3}\right)^{\!\frac{1}{3}}$ in the case of the Whitham equation).

Equation~\eref{eq:steadyproblem} also admits a variational formulation: local minimisers of the functional
$\Efunc: H^1(\R) \rightarrow \R$ given by
\begin{equation}\label{eq:Edefn}
\Efunc(u) = \underbrace{-\frac{1}{2}\int_\R uLu \dx}_{\displaystyle :=\Lfunc(u)}\, \underbrace{- \int_\R N(u) \dx}_{\displaystyle :=\Nfunc(u)},
\end{equation}
where $N$ is the primitive function of $n$ which vanishes at the origin, so that
$$
N(x) := N_{p+1}(x)+N_\mathrm{r}(x), \qquad N_\mathrm{r}(x):=\int_0^x n_\mathrm{r}(s) \ds,
$$
under the constraint that $\Qfunc$ is held fixed are solitary-wave solutions of~\eref{eq:steadyproblem}.
The technique employed by Albert and Zeng, which relies upon the fact that $L$ is of positive order
(so that $\Efunc_\mathrm{lw}$ is coercive), is however not applicable in the present situation in which $L$ is a smoothing operator. Instead we use methods developed by Buffoni \cite{Buffoni04a} and 
Groves \& Wahl\'{e}n~\cite{GrovesWahlen11}. We consider a fixed ball
$$U=\{u \in H^1(\R): \|u\|_1 < R\},$$
and seek small-amplitude solutions, that is solutions in the set
$$
U_{\mu} := \left\{ u \in U:\; \Qfunc(u) = \mu \right\},
$$
where $\mu$ is a small, positive, real number. In particular we examine minimising sequences for $\Efunc$ over $U_\mu$
which do not approach the boundary of $U$, and establish the following result with the help
of the concentration-compactness principle.

\begin{theorem}[Existence]\label{theorem:main}
There exists $\mu_\star>0$ such that the following statements hold for each $\mu \in (0,\mu_\star)$.
\begin{itemize}
\item[(i)] The set $D_\mu$ of minimisers of $\Efunc$ over the set $U_\mu$
is non-empty and the estimate $\|u\|_1^2 = \bigO(\mu)$ holds uniformly over $u \in D_\mu$
and $\mu \in (0,\mu_\star)$.
Each element of $D_\mu$ is a solution of the travelling-wave equation~\eref{eq:steadyproblem};
the wave speed $\nu$ is the Lagrange multiplier in this constrained variational principle.
The corresponding solitary waves are supercritical, that is their speed $\nu$ exceeds $m(0)$.\\[-6pt]
\item[(ii)] Let $s < 1$ and suppose that $\{u_n\}_{n\in\N_0}$ is a minimising sequence for $\Efunc$ over
$U_\mu$ with the property that 
\begin{equation}
\label{eq:away from boundary}
\sup_{n \in \N_0}\|u_n\|_1 < R.
\end{equation}
There exists a sequence $\{x_n\}_{n\in \N_0}$ of real numbers such that a subsequence of $\{u_n(\cdot + x_n)\}_{n\in\N_0}$ converges in $H^s(\R)$ to a function in $D_\mu$.
\end{itemize}
\end{theorem}

Theorem~\ref{theorem:main} is proved in two steps. We begin by constructing a minimising sequence which
satisfies condition~\eref{eq:away from boundary}. To this end
we consider the corresponding problem for periodic travelling waves (see Section~\ref{sec:periodic})
and penalise the variational functional so that minimising sequences do not approach the boundary of the corresponding domain
in function space. Standard methods from the calculus of variations yield the existence of minimisers for the
penalised problem, and \emph{a priori} estimates confirm that the minimisers lie in the region unaffected by the
penalisation; in particular they are bounded (uniformly over all large periods) away from the boundary.
A minimising sequence $\{\tilde{u}_n\}_{n \in \N_0}$ for $\Efunc$ over $U_\mu$ is obtained by letting the period tend to infinity.

The minimising sequence $\{\tilde{u}_n\}_{n \in \N_0}$ is used to show that the quantity
$$I_\mu :=\inf \left\{ \Efunc(u) \colon u \in U_{\mu} \right\}$$
is \emph{strictly subadditive}, that is
$$
I_{\mu_1 + \mu_2} < I_{\mu_1} + I_{\mu_2} \quad \mbox{ whenever }0 < \mu_1,\mu_2  < \mu_1 + \mu_2 < \mu_\star.
$$
The proof of this fact, which is presented in Section~\ref{sec:SSA},
is accomplished by showing that the functions $\tilde{u}_n$ `scale'
in a fashion similar to the long-wave \emph{Ansatz}~\eref{eq:longwaveansatz}; we may therefore approximate $\Efunc$
by a scaling of $\Efunc_\mathrm{lw}$ along this minimising sequence. The corresponding strict subadditivity result for the
latter functional is a straightforward matter, and a perturbation argument shows that it remains valid for $\Efunc$.

In a second step we apply the concentration-compactness principle to show that \emph{any} minimising sequence
satisfying~\eref{eq:away from boundary} converges --- up to subsequences and translations --- in
$H^s(\R)$, $s<1$ to a minimiser of $\Efunc$ over $U_\mu$ (Section~\ref{sec:concentration}). The strict subadditivity
of $I_\mu$ is a key ingredient here. The proof of Theorem~\ref{theorem:main}(i) is completed by \emph{a priori}
estimates for the size and speed of solitary waves obtained in this fashion.

Section~\ref{sec:consequences} examines some consequences of Theorem~\ref{theorem:main}. In particular, the
relationship between the solutions to~\eref{eq:lwe} found in Theorem~\ref{theorem:AlbertZeng} and the
solutions to~\eref{eq:steadyproblem} found in Theorem~\ref{theorem:main} is rigorously clarified.
Under an additional regularity hypothesis upon $n$ we show that every solution $u$ in the
set $D_\mu$ lies in $H^{2j_\star}(\R)$,
`scales' according to the long-wave \emph{Ansatz}~\eref{eq:longwaveansatz} and satisfies
$$
\mathrm{dist}_{H^{j_\star}(\R)} \big(\mu^{-\alpha} u( \mu^{-\beta}\cdot) , D_\mathrm{lw} \big) \to 0
$$
as $\mu \searrow 0$; the convergence is uniform over $D_\mu$. Corresponding convergence results for the wave speeds
and infima of $\Efunc$ over $U_\mu$ and $\Efunc_\mathrm{lw}$ over $\{w \in H^{j_\star}(\R): \Qfunc(w)=1\}$ are
also presented. These results may contribute towards
the discussion of the validity of the Whitham equation as a model for water waves: they show that the Whitham solitary waves
are approximated by Korteweg-deVries solitary waves, and it is known that solutions of the Korteweg-deVries equation
do approximate the solutions of the full water-wave problem (Craig~\cite{Craig85},
Schneider \& Wayne~\cite{SchneiderWayne00}).   

Theorem~\ref{theorem:main} also yields information about the stability of the set of solitary-wave solutions
to~\eref{eq:DEproblem} defined by $D_\mu$. Observing that $\Efunc$ and $\Qfunc$ are conserved
quantities associated with equation~\eref{eq:DEproblem}, we apply a general principle that the solution
set of a constrained minimisation problem of this type constitutes a stable set of solutions of the corresponding
initial problem (Theorem \ref{theorem:stability}): choosing $\mathrm{dist}_{L^2(\R)}(u(0), D_\mu)$ sufficiently small
ensures that $\mathrm{dist}_{L^2(\R)}(u(t),D_\mu)$ remains small over the time of existence of a solution
$u: [0,T] \rightarrow H^1(\R)$ with $\sup_{t \in [0,T]} \|u(t)\|_1 < R$. Of course the well-posedness
of the initial-value problem for equation~\eref{eq:DEproblem} is a prerequisite for discussing the stability of $D_\mu$.
This discussion is however outside the scope of the present paper; we merely assume that the initial-value problem
is locally well posed in a sense made precise in Section~\ref{sec:consequences}. 
Our stability result is \emph{conditional} since it
applies to solutions only for as long as they remain in $U$ (for example certain solutions of the
Whitham equation~\eref {eq:whitham} have only a finite time of existence
(Naumkin \& Shishmarev~\cite{NaumkinShishmarev})), and \emph{energetic} since distance
is measured in $L^2(\R)$ rather than $H^1(\R)$ (note that the norms in $H^s(\R)$ for $s\in[0,1)$ are
all metrically equivalent on $U$). Theorem~\ref{theorem:stability} also refers to the stability of the entire set
$D_\mu$; in the special case where the minimiser of $\Efunc$ over $U_\mu$ is unique up to translations
it coincides with (conditional and energetic) orbital stability of this solution.




\section{Preliminaries}\label{sec:preliminaries}   

\subsection*{Functional-analytic setting for the solitary-wave problem}

Let $\SS(\R)$ be the Schwartz space of rapidly decaying smooth functions, and let
$\FF$ denote the unitary Fourier transform on $\SS(\R)$, so that
$$\FF(\varphi)(k) := \frac{1}{\sqrt{2\pi}} \int_\R \varphi(x) \exp(-\mathrm{i}kx)\dx,  \qquad \varphi \in \SS(\R),$$
and on the dual space of tempered distributions $\SS^\prime(\R)$, so that
$(\hat{f},\varphi)=(f,\hat{\varphi})$ for $f \in \SS^\prime(\R)$. 
By $L^p(\R)$, $p \geq 1$ we denote the space of real-valued $p$-integrable functions with norm
$\| f \|_{L^p(\R)} : = \int_\R |f(x)|^p \dx$,
by $H^s(\R)$, $s \in \R$ the real Sobolev space consisting of those tempered distributions for which
the norm
$$
\| f \|_s := \left( \int_\R |\hat f(k)|^2 \left( 1 + k^2 \right)^{s}\dk \right)^{\!\!\frac{1}{2}}
$$
is finite, and by $\mathrm{BC}(\R)$ the space of bounded and continuous real-valued functions with finite supremum norm $\|f\|_\infty := \sup_{x \in \R} |f(x)|$; there is a continuous embedding $H^s(\R) \hookrightarrow \mathrm{BC}(\R)$ for any $s > \frac{1}{2}$, so that
$\|u\|_\infty \leq c_s\|u\|_s$ for all $u \in H^s(\R)$.
We write $L^2(\R)$ for $H^0(\R)$, and for all spaces the subscript `c' denotes the subspace of compactly supported functions, so that
$$
H_\mathrm{c}^s(\R) := \{ f \in H^s(\R) \colon \mbox{$\supp(f)$ is compact} \}.
$$

We now list some basic properties of the operators $L$, $n$ appearing in equation~\eref{eq:steadyproblem}
and functionals $\Efunc$, $\Qfunc$ defined in equations~\eref{eq:Qdefn}, \eref{eq:Edefn}.

\begin{proposition}\label{prop:Ln1} \hspace{1in}
\begin{itemize}
\item[(i)] The linear operator $L$ belongs to $C^\infty(H^s(\R),H^{s+|m_0|}(\R)) \cap C^\infty(\SS(\R),\SS(\R))$ for each $s \geq 0$.\\[-10pt] 
\item[(ii)] For each $j \in \N$ there exists a constant $\tilde C_j > 0$ such that
$$
| Lu(x) | \leq \frac{\tilde C_j}{\left(\dist (x, \supp(u) ) \right)^j} \, \| u\|_0, \qquad x \in \R \setminus \supp(u)
$$
for all $u \in L_\mathrm{c}^2(\R)$.
\item[(iii)] Suppose that $n \in C^k(\R, \R)$ for some $k \in \N$. For each $R>0$ the
function $n$ induces a continuous Nemitskii operator $B_R(0) \subset H^s(\R) \rightarrow H^s(\R)$, where $s \in (\frac{1}{2},k]$.\\[-10pt]
\end{itemize}
\end{proposition}
\proof (i) This assertion follows directly from the definition of $L$.

(ii) Assumption~\eref{eq:mclassical} implies that $m^{(j)} \in L^2(\R)$ for any ${j \in \N}$. Applying Plancherel's theorem and H\"older's inequality to the convolution formula~\eref{eq:Lconvolution}, one finds that
\begin{eqnarray*}
| Lu(x) | &=& \frac{1}{\sqrt{2\pi}} \left| \int_{\supp(u)} \frac{(x-y)^j}{(x-y)^j} K(x-y) u(y)\,dy \right| \\
& \leq &\frac{C_j  \| m^{(j)}\|_0}{\sqrt{2\pi}} \, \frac{1}{\dist(x,\supp(u))^j}\,  \|u\|_0.
\end{eqnarray*}

(iii) Construct a $k$ times continuously
differentiable function $\tilde{n}: \R \rightarrow \R$ whose derivatives are bounded and which satisfies
$\tilde{n}(x) = n(x)$ for $|x| \leq c_s R$ (for example by multiplying
$n$ by a smooth `cut-off' function). The results given by Bourdaud \& Sickel~\cite[Theorem 7]{BourdaudSickel11}
(for $s \in (0,1)$) and Brezis \& Mironescu~\cite[Theorem 1.1]{BrezisMironescu01} (for $s \geq 1$) show that
$\tilde{n}$ induces a continuous Nemitskii operator $H^s(\R) \to H^s(\R)$ for $s \in (0, k]$ and hence that $n$ induces a continuous Nemitskii operator $B_R(0) \subset H^s(\R) \rightarrow H^s(\R)$ for $s \in (\frac{1}{2}, k]$.\qed

According to the previous proposition we may study~\eref{eq:steadyproblem} as an equation in $H^s(\R)$
for $s>\frac{1}{2}$ (provided that $n$ is sufficiently regular). In keeping with this observation we work
in the fixed ball $U=\{u \in H^1(\R): \|u\|_1 < R\}$.

\begin{proposition}\label{prop:Ln2}
Suppose that $n \in C^1(\R)$.
\begin{itemize}
\item[(i)]
The functionals $\Lfunc$, $\Nfunc$ and $\Qfunc$ belong to $C^1(U,\R)$
and their $L^2(\R)$-derivatives are given by the formulae
$$
\Lfunc^\prime(u) :=  -Lu, \quad \Nfunc^\prime(u) := -n(u), \quad \Qfunc^\prime(u) = u.
$$
These formulae define functions $\Lfunc^\prime$, $\Nfunc^\prime$, $\Qfunc^\prime \in C(U,H^1(\R))$.
\item[(ii)]
The functional $\Efunc$ belongs to $C(H^s(\R),\R)$ for each $s > \frac{1}{2}$.
\end{itemize}
\end{proposition}\smallskip

Finally, we note that solutions of the travelling-wave equation may inherit further regularity from $n$.

\begin{lemma}[Regularity] \label{lemma:regularity}
Suppose that $n \in C^{k+1}(\R)$ for some $k \in \N$. For sufficiently small values of $R$, 
every solution $u \in U$ of~\eref{eq:steadyproblem} belongs to $H^{k+1}(\R)$ and
satisfies
$$\|u\|_{k+1} \leq c \|u\|_1.$$
\end{lemma}
\proof
Differentiating~\eref{eq:steadyproblem}, we find that 
\begin{equation}\label{eq:diffproblem}
u^\prime = \frac{L u^\prime}{\nu- n^\prime(u)}.
\end{equation}
There exists a positive constant $c_\delta$ such that $\nu - n^\prime(u) \geq \delta > 0$ whenever $\|u\|_\infty < c_\delta$;
the embedding
$H^1(\R) \hookrightarrow \mathrm{BC}(\R)$ guarantees that this condition is fulfilled for each $u \in U$
for sufficiently small values of $R$.

Suppose that $m \in \{1,\ldots, k\}$. For each fixed $u \in H^m(\R)$ the formula
$$
\varphi_u(v) = \frac{v}{\nu- n^\prime(u)}
$$
defines an operator in $B(L^2(\R), L^2(\R))$ and $B(H^m(\R),H^m(\R))$, and by interpolation it
follows that $\varphi_u \in B(H^s(\R),H^s(\R))$ for $s \in [0,m]$; its norm depends upon $\|u\|_m$.
Furthermore, recall that $L \in B(H^s(\R),H^{s+|m_0|}(\R))$ for all $s \in [0,\infty)$, so that
$$
\psi_u:=\varphi_u \circ L \in B\left(H^s(\R), H^{s_\star}(\R)\right), \qquad s_\star=\min(m,s+|m_0|),
$$
and the norm of $\psi_u$ depends upon $\|u\|_m$.

It follows that any solution $w \in H^s({\mathbb R})$ of the equation
\begin{equation}\label{eq:bs}
w = \psi_u (w)
\end{equation}
in fact belongs to $H^{s_\star}(\R)$, where $s_\star = \min(m,s+|m_0|)$, and satisfies the estimate
$$\|w\|_{s_\star} \leq c_{\|u\|_m} \|w\|_s.$$
Applying this argument recursively, one finds that any solution $w \in L^2(\R)$ of~\eref{eq:bs} belongs to
$H^m(\R)$ and satisfies
$$\|w\|_m \leq c_{\|u\|_m} \|w\|_0.$$

Observe that equation~\eref{eq:diffproblem} is equivalent to $u^\prime = \psi_u (u^\prime)$. A bootstrap
argument therefore shows that $u^\prime \in H^k(\R)$ with
$$\|u^\prime\|_m \leq c_{\|u\|_1} \|u^\prime\|_0, \qquad m=1,\ldots,k.\eqno{\Box}$$

\subsection*{Functional-analytic setting for the periodic problem}

Let $P > 0$. Let $L_P^2$ be the space of $P$-periodic, locally square-integrable functions with
Fourier-series representation
$$
u(x) = \frac{1}{\sqrt{P}} \sum_{k \in \Z} \hat u_k \exp(2\pi \mathrm{i} kx/P),
$$  
and define
$$
H_P^s := \left\{ u \in L_P^2 \colon \|u \|_{H_P^s} :=  \left( \sum_{k \in \Z} \left( 1 + \frac{4\pi^2 k^2}{P^2}\right)^{\!\!s} | \hat u_k |^2 \right)^{\!\!\frac{1}{2}}  < \infty \right\}
$$
for $s \geq 0$. Just as for the Sobolev spaces $H^s(\R)$ one has the continuous embedding
$H_P^s \hookrightarrow \mathrm{BC}(\R)$ for all $s > \frac{1}{2}$; the embedding constant is
independent of $P$.


\begin{proposition}
The operator $L$ extends to an operator $\SS^\prime(\R) \to \SS^\prime(\R)$ which maps $H_P^s$ smoothly into $H_P^{s+|m_0|}$, acting on the Fourier coefficients $\hat{u}_k$, $k \in \Z$,
of a function $u$ by pointwise multiplication, so that
$$
(\widehat{Lu})_k = m (2\pi k / P)\, \hat u_k, \qquad k \in \Z.
$$
\end{proposition}
\proof
The operator $L$ is symmetric on $L^2(\R)$ and maps $\SS(\R)$ into itself; it therefore extends to an operator $\SS^\prime(\R) \to \SS^\prime(\R)$. In particular, the convolution theorem shows that $L$ maps $P$-periodic functions to $P$-periodic functions, acting on their Fourier coefficients by pointwise multiplication;
it follows that $L \in C^\infty(H_P^s, H_P^{s+|m_0|})$.\qed

There is a natural injection from the set of functions $\tilde u_P \in L_\mathrm{c}^2(\R)$ with $\supp (\tilde u) \subset (-\frac{P}{2},\frac{P}{2})$
to $L_P^2$, namely 
$$
\tilde u_P \mapsto u_P := \sum_{j \in \Z} \tilde u_P (\cdot + jP),
$$
where the series converges in $\SS^\prime(\R)\, \cap\, L^2_\mathrm{loc}(\R)$.
The following proposition shows that this map commutes with $L$.

\begin{proposition}
Any function $\tilde u_P \in L_\mathrm{c}^2(\R)$ with $\supp (\tilde u_P) \subset (-\frac{P}{2},\frac{P}{2})$ satisfies
$$
L^2_\mathrm{c}
 \ni \sum_{|j| \leq J} L \tilde u_P (\cdot + jP) \stackrel{J \to \infty}{\longrightarrow} L u_P \in L_P^2
 $$
in $\SS^\prime(\R)\, \cap\, L^2_\mathrm{loc}(\R)$.
\end{proposition}
\proof The convergence in $\SS^\prime(\R)$ follows from the continuity of
$L: \SS^\prime(\R) \to \SS^\prime(\R)$, while that in $L_\mathrm{loc}^2(\R)$
follows from the calculation
\begin{eqnarray*}
\lefteqn{\Big\| \sum_{|j| \geq J} L\tilde u_P(\cdot + jP) \Big\|_{L^2((-M,M))}}\qquad\\ 
& \leq & \sum_{|j| \geq J} \| L\tilde u_P(\cdot + jP) \|_{L^2((-M,M))}\\
& \leq & (2M)^\frac{1}{2} \tilde C_2 \|\tilde u_P \|_0 \sum_{|j| \geq J} \frac{1}{\dist\left([-M,M], \supp (\tilde u (\cdot + jP)) \right)^{2}}\\ 
& \leq &(2M)^\frac{1}{2} \tilde C_2 \|\tilde u_P\|_0 \sum_{|j| \geq J} \frac{1}{((|j|-\frac{1}{2})P - M)^2}\\
& \to & 0
\end{eqnarray*}
as $J \to \infty$.\qed

Define
$$
U_P := \left\{ u \in H_P^1 \colon \|u\|_{H_P^1} < R \right\}
$$
and functionals $\Nfunc_P$, $\Lfunc_P$, $\Efunc_P$, $\Qfunc_P \colon U_{P} \to \R$ by replacing the domain of integration in the definitions of $\Nfunc$, $\Lfunc$, $\Efunc$, $\Qfunc$ by one period $(-\frac{P}{2},\frac{P}{2})$. Observing that Proposition~\ref{prop:Ln2} (with the obvious modifications) holds for the new functionals,
we study $\Efunc_P$, $\Qfunc_P \in C^1(U_P,\R)$. Each minimiser of $\Efunc_P$ over the set
$$
U_{P,\mu} := \left\{ u \in U_{P} \colon \Qfunc_P(u) = \mu \right\}
$$
is a $P$-periodic solution of the travelling-wave equation~\eref{eq:steadyproblem}; 
the wave speed $\nu$ is the Lagrange multiplier in this constrained variational principle.

\subsection*{Additional notation}

\begin{itemize}
\item
We denote the set of functions which are square integrable over an open subset $S$ of ${\mathbb R}$
by $L^2(S)$ and the subset of $L^2(S)$ consisting of those functions whose weak derivative
exists and is square integrable by $H^1(S)$.
\item
The symbol $c$ denotes a a generic constant which is independent of $\mu \in (0,\mu_\star)$ (and of
course functions in a given set or sequence); its dependence
upon other quantities is indicated by a subscript. All order-of-magnitude estimates are also
uniform over $\mu \in (0,\mu_\star)$, and in general we replace $\mu^\star$ with a smaller number
if necessary for the validity of our results.
\end{itemize}




\section{The minimisation problem for periodic functions}\label{sec:periodic}

\subsection*{The penalisation argument}

Seeking a constrained minimiser of $\Efunc_P$ in the set
$U_{P,\mu}$ by the direct method of the
calculus of variations, one is confronted by the difficulty that a minimising sequence
may approach the boundary of $U_P$. To overcome this difficulty we observe
that $\Efunc_P$ also defines a continuously differentiable functional
on the set
$$V_P:=\{ u \in H_P^1 \colon \|u\|_{H_P^1} < 2R \}$$
and consider the auxiliary functional
$$
\Efunc_{P,\varrho}(u) := \Efunc_P(u) + \varrho\left(\|u\|_{H_P^1}^2\right)
$$
with constraint set
$$
V_{P,\mu} := \left\{ u \in H_P^1 \colon \|u\|_{H_P^1} < 2R, \; \Qfunc(u) = \mu  \right\},
$$
where we note the helpful estimate
\begin{equation}\label{eq:helpfulest}
\| u \|_\infty \leq c \| u \|_{L_P^2}^\frac{1}{2} \|  u \|_{H_P^1}^\frac{1}{2} \leq  c \mu^\frac{1}{4}, \qquad u \in V_{P,\mu}.
\end{equation}

Here $\varrho\colon [0,(2R)^2) \to [0,\infty)$ is a smooth, increasing `penalisation' function such that\\[-8pt]
\begin{itemize}
\item[(i)] $\varrho(t) = 0$  whenever  $0 \leq t \leq R^2$,\\[-8pt]
\item[(ii)] $\varrho(t) \to \infty$  as  $t \nearrow (2R)^2$,\\[-8pt]
\item[(iii)] for every constant $a_1 \in (0,1)$ there exist $M_1, M_2 > 0$ and $a_2 > 1$ such that 
$$
\varrho^\prime(t) \leq M_1 (\varrho(t))^{a_1} + M_2 (\varrho(t))^{a_2};
$$
\end{itemize}
an example of such a function $\varrho$ can be obtained by scaling and translating the function
$$
t \mapsto 
\left\{\begin{array}{ll}
(1-t)^{-1} \exp(-1/t), \qquad &t \in (0,1),\\[1mm]
0, \qquad & t \leq 0.
\end{array}\right.
$$


The following lemma is obtained by standard weak continuity arguments (e.g.\ see Struwe~\cite[\S\S I.1, I.2]{Struwe}). 

\begin{lemma}
The functional $\Efunc_{P,\varrho} \colon V_{P,\mu} \to \R$ is weakly lower semicontinuous, bounded from below, and satisfies
$\Efunc_{P,\varrho}(u) \to \infty$ as $\|u\|_{H_P^1} \nearrow 2R$. In particular, it
has a minimiser $\bar u_{P} \in V_{P,\mu}$.
\end{lemma}

The next step is to show that $\bar u_P$ in fact minimises
$\Efunc_P$ over $U_{P,\mu}$. This result
relies upon estimates for $\Efunc_{P,\varrho}$ which
are uniform in $P$ and are derived in Lemmata~\ref{lemma:umu}
and~\ref{lemma:PE-convergence} and Corollary~\ref{cor:umu}
by examining the functional $\Efunc$ and its relationship to $\Efunc_{P,\varrho}$.


\begin{lemma}\label{lemma:umu}
For any $w \in W$ the `long-wave test function' $S_\mathrm{lw}w$, where
$$(S_\mathrm{lw}w)(x) = \mu^\alpha w(\mu^\beta x),$$
lies in $U$ and satisfies
$$
\Efunc(S_\mathrm{lw}w) = -\mu m(0) + \mu^{1+(p-1)\alpha} \Efunc_\mathrm{lw}(w) + o(\mu^{1+(p-1)\alpha}),
$$
where the values of $\alpha$ and $\beta$ are given by~\eref{eq:alphabeta} and
$\Efunc_\mathrm{lw}$ is defined in equation~\eref{eq:redEdefn}. The estimate holds uniformly
over $w \in W$, and $w \in W_1$ implies $u \in U_\mu$.
\end{lemma}
\proof
Observe that
$$\Qfunc(S_\mathrm{lw}w) = \mu^{2\alpha-\beta}, \qquad \FF[S_\mathrm{lw}w](k) = \mu^{\alpha-\beta} \hat w(\mu^{-\beta} k)$$
and
$$\|S_\mathrm{lw}w\|_1^2 = \mu^{2\alpha -\beta}\|w\|_0^2 + \mu^{2\alpha + \beta}\|w^\prime\|_0^2 \leq c \mu$$
for $\alpha, \beta>0$ with $2\alpha - \beta \geq 1$. A direct calculation shows that
\begin{eqnarray*}
\afl \Efunc(S_\mathrm{lw}w) & = & - \frac{1}{2} \int_\R m(k) \, |\FF[S_\mathrm{lw}w](k)|^2 \dk-  \mu^{-\beta} \int_\R N\left(\mu^{\alpha} w(x)\right)\dx \\ 
\afl &=  & -\mu^{2\alpha - \beta} m(0) -   \frac{\mu^{2\alpha+(2j_\star - 1)\beta}\, m^{(2j_\star)}(0)}{2 (2j_\star)!}   \int_{\R}  k^{2j_\star} |\hat w(k)|^2 \dk\\
\afl &  &  \quad \mbox{}- \mu^{(p+1)\alpha-\beta} \int_\R N_{p+1}(w(x)) \dx \\
\afl & & \quad \mbox{} - \mu^{-\beta} \int_\R N_\mathrm{r}(\mu^\alpha w(x)) \dx
- \frac{\mu^{2\alpha-\beta}}{2} \int_\R r(\mu^\beta k) |\hat w(k)|^2 \dk,
\end{eqnarray*}
and one can estimate
\begin{eqnarray*}
\afl\left| \frac{\mu^{2\alpha-\beta}}{2} \int_\R r(\mu^\beta k) |\hat w(k)|^2 \dk + \mu^{-\beta} \int_\R N_\mathrm{r}(\mu^\alpha w(x)) \dx\right| \\
\afl\qquad\leq c\left(\mu^{2\alpha+ (2j_\star+1)\beta} \int_\R  k^{2 j_\star +2} |\hat w(k)|^2 \dk + \mu^{(p+\delta+1)\alpha-\beta} \int_\R |w(x)|^{p+\delta+1} \dx \right).
\end{eqnarray*}
Choosing $\alpha$ and $\beta$ such that $(p-1)\alpha= 2j_\star\beta$ and $2\alpha - \beta = 1$, so that
$\alpha$ and $\beta$ are given by~\eref{eq:alphabeta},
yields the desired estimate.\qed


\begin{lemma}\label{lemma:PE-convergence}
Let $\{\tilde u_P\}_P$ be a bounded family of functions in $H^1(\R)$ with 
$$
\supp(\tilde u_P) \subset (-\tfrac{P}{2},\tfrac{P}{2}) \quad\mbox{ and }\quad\dist(\pm \tfrac{P}{2},\supp(\tilde u_P)) \geq \frac{1}{2} P^\frac{1}{4}
$$
and define $u_P \in H_1^P$ by the formula
$$
u_P = \sum_{j\in \Z} \tilde u_P(\cdot + jP).
$$
\begin{itemize}
\item[(i)]
The function $u_P$ satisfies
$$
\lim_{P \to \infty} \| L \tilde u_P - L u_P \|_{H^1(-\frac{P}{2},\frac{P}{2})} = 0, \qquad
\lim_{P \to \infty} \|L\tilde u_P \|_{H^1(\{|x| > \frac{P}{2}\})} = 0,
$$
\item[(ii)]
The functionals $\Efunc$, $\Qfunc$ and $\Efunc_P$, $\Qfunc_P$ have the properties that
$$
\lim_{P \to \infty} \left( \Efunc(\tilde u_P) - \Efunc_P(u_P) \right) = 0, \qquad \Qfunc(\tilde u_P) = \Qfunc_P(u_P)
$$
and
$$
\lim_{P \to \infty} \| \Efunc^\prime (\tilde u_P) - \Efunc^\prime_P (u_P) \|_{H^1(-\frac{P}{2},\frac{P}{2})} = 0, \qquad
\lim_{P \to \infty} \|\Efunc^\prime (\tilde u_P) \|_{H^1(\{|x| > \frac{P}{2}\})}=0,
$$
$$
\| \Qfunc^\prime (\tilde u_P) - \Qfunc^\prime_P (u_P) \|_{H^1(-\frac{P}{2},\frac{P}{2})} = 0, \qquad
\|\Qfunc^\prime (\tilde u_P) \|_{H^1(\{|x| > \frac{P}{2}\})}=0.
$$
\end{itemize}
\end{lemma} 
\proof
(i) Using Proposition~\ref{prop:Ln1}(ii), we find that
\begin{eqnarray*}
\lefteqn{\int_{-\frac{P}{2}}^{\frac{P}{2}} | L\tilde u_P - L u_P|^2\dx}\qquad\\ 
&=& \int_{-\frac{P}{2}}^{\frac{P}{2}} \Big| \sum_{|j| \geq 1} L \tilde u_P(x+jP) \Big|^2\dx\\
&\leq &\int_{-\frac{P}{2}}^{\frac{P}{2}} \left( \sum_{|j| \geq 1} \frac{\tilde C_3 \|\tilde u_P\|_0}{\dist\left(x+jP, \supp(\tilde u_P) \right)^3} \right)^{\!\!2}\dx\\
&\leq &\int_{-\frac{P}{2}}^{\frac{P}{2}} \left( 2 \sum_{j \geq 0} \frac{\tilde C_3 \|\tilde u_P\|_0}{\left( jP +  \frac{1}{2} P^\frac{1}{4} \right)^3} \right)^{\!\!2}\dx\\
&\to & 0,\\
\\
\lefteqn{\int_{|x| > \frac{P}{2}} |L\tilde u_P|^2 \dx}\qquad\\ 
&\leq & \tilde C_1^2 \|\tilde u_P\|_0^2 \int_{|x| > \frac{P}{2}} \frac{\dx}{\dist\left( x, \supp( \tilde u_P) \right)^2} \\
&\leq & \tilde C_1^2 \|\tilde u_P\|_0^2 \int_{|x| > \frac{P}{2}} \frac{\dx}{\left( |x| - \frac{1}{2}(P - P^\frac{1}{4})\right)^2}\\ 
&= & \frac{4 \tilde C_1^2 \|\tilde u_P\|_0^2}{P^\frac{1}{4}} \\
& \to & 0
\end{eqnarray*}
and therefore
$$
\lim_{P \to \infty} \| L \tilde u_P - L u_P \|_{L^2(-\frac{P}{2},\frac{P}{2})} = 0, \qquad
\lim_{P \to \infty} \|L\tilde u_P \|_{L^2(\{|x| > \frac{P}{2}\})} = 0,
$$
as $P \to \infty$. The same calculation is valid with $u_P$ and $\tilde u_P$ replaced by respectively
$u_P^\prime$ and $\tilde{u}_P^\prime$, and since $L$ commutes with differentiation this observation completes
the proof.

(ii) Observe that
\begin{eqnarray*}
\lefteqn{|\Lfunc(\tilde u_P) - \Lfunc_P(u_P)|} \qquad\\
& = & \left| \frac{1}{2}\int_\R \tilde u_P L \tilde u_P \dx - \frac{1}{2}\int_{-\frac{P}{2}}^{\frac{P}{2}} u_P L u_P\dx \right| \\
&=& \left|\frac{1}{2} \int_{-\frac{P}{2}}^{\frac{P}{2}} \tilde u_P (L\tilde u_P - L u_P)\dx \right| \\
&\leq & \frac{1}{2}\| \tilde u_P \|_0 \|L\tilde u_P - L u_P \|_{L^2(-\frac{P}{2},\frac{P}{2})} \\
& \to & 0
\end{eqnarray*}
and
$$
 \| \Lfunc^\prime (\tilde u_P) - \Lfunc^\prime_P (u_P) \|_{H^1(-\frac{P}{2},\frac{P}{2})} = \| L \tilde u_P - L u_P \|_{H^1(-\frac{P}{2},\frac{P}{2})} \to 0,
$$
$$
\|\Lfunc^\prime (\tilde u_P) \|_{H^1(\{|x| > \frac{P}{2}\})}=\|L\tilde u_P \|_{H^1(\{|x| > \frac{P}{2}\})} \to 0
$$
as $P \rightarrow \infty$. Furthermore
$$
\Nfunc(\tilde u_P) = -\int_\R N(\tilde u_P)\dx = -\int_{-\frac{P}{2}}^{\frac{P}{2}} N(\tilde u_P)\dx =  -\int_{-\frac{P}{2}}^{\frac{P}{2}} N(u_P)\dx = \Nfunc_P(u_P)
$$
and
\begin{eqnarray*}
\afl & & \Nfunc^\prime(\tilde u_P(x)) = -n(\tilde u_P(x)) =
\left\{\begin{array}{ll} -n(u_P(x))=\Nfunc_P^\prime(u_P(x)), & x \in (-\frac{P}{2},\frac{P}{2}),\\[1mm]
0, & |x| \geq \frac{P}{2},
\end{array}\right.\\[3mm]
\afl & & (\Nfunc^\prime(\tilde u_P))^\prime(x) = -n^\prime(\tilde u_P(x))\tilde{u}_P^\prime(x) \\[1mm]
\afl & & \hphantom{(\Nfunc^\prime(\tilde u_P))^\prime(x) }=
\left\{\begin{array}{ll} -n^\prime(u_P(x))u_P^\prime(x)=(\Nfunc_P^\prime(u_P))^\prime(x), & x \in (-\frac{P}{2},\frac{P}{2}),\\[1mm]
0, & |x| \geq \frac{P}{2},
\end{array}\right.
\end{eqnarray*}
so that
$$
\| \Nfunc^\prime (\tilde u_P) - \Nfunc^\prime_P (u_P) \|_{H^1(-\frac{P}{2},\frac{P}{2})} = 0, \qquad
\|\Nfunc^\prime (\tilde u_P) \|_{H^1(\{|x| > \frac{P}{2}\})}=0.
$$
The result for $\Efunc$, $\Efunc_P$ follows from these calculations and the formulae
$\Efunc = \Lfunc+\Nfunc$, $\Efunc_P = \Lfunc_P+\Nfunc_P$, and a similar calculation yields the result for
$\Qfunc$, $\Qfunc_P$.\qed


\begin{corollary}\label{cor:umu}
There exist constants $I_\star > 0$ and $P_{\mu}>0$ such that
$$I_\mu:=\inf\left\{\Efunc(u) \colon u \in U_{\mu} \right\} < - \mu m(0) - \mu^{1+(p-1)\alpha} I_\star$$
and
$$
I_{P,\varrho,\mu} := \inf \left\{ \Efunc_{P,\varrho}(u) \colon u \in V_{P,\mu} \right\} < - \mu m(0) -  \mu^{1+(p-1)\alpha} I_\star
$$
for each $P \geq P_{\mu}$. 
\end{corollary}
\proof
Taking $\psi \in C_\mathrm{c}^\infty(\R)$ with $\Qfunc(\psi)=1$ and writing $w(x)=\sqrt{\lambda} \psi(\lambda x)$,
one finds that
$$\Efunc_\mathrm{lw}(w) = -\lambda^{2j_\star}\frac{m^{(2j_\star)}(0)}{2(2j_\star)!}\int_\R (\psi^{(j_\star)})^2\dx-\lambda^{(p-1)/2}\int_\R N_{p+1}(\psi)\dx<0$$
for sufficiently small values of $\lambda$ provided that $p<4j_\star+1$
and $N_{p+1}(\psi)>0$; these conditions are satisfied under
assumption (A3) by choosing $\psi>0$ if $c_p>0$ and $\psi<0$ if $c_p<0$.

Noting that $w \in W$ for sufficiently large values of $S$, we find from Lemma~\ref{lemma:umu} that
\begin{eqnarray}\label{eq:Itildestar}
\Efunc(S_\mathrm{lw}w) + \mu m(0) & = & 
\mu^{1+(p-1)\alpha} \Efunc_\mathrm{lw}(w) + o(\mu^{1+(p-1)\alpha}) \nonumber \\
& < & \tfrac{1}{2}\mu^{1+(p-1)\alpha} \Efunc_\mathrm{lw}(w).
\end{eqnarray}
Observe that $\supp(S_\mathrm{lw}w) = \mu^{-\beta} \supp(w)$, so that $S_\mathrm{lw}w$ satisfies the assumptions of Lemma~\ref{lemma:PE-convergence} if and only if $\mu^\beta P \geq c_w$, where $c_w$ is a positive constant independent of $P$.    
For such $P$ a combination of Lemma~\ref{lemma:PE-convergence} and~\eref{eq:Itildestar} yields
\begin{eqnarray*}
I_{P, \varrho, \mu} &\leq& \Efunc_P(u_P)\\ 
&\leq & -\mu m(0) + \tfrac{1}{2}\mu^{1+(p-1)\alpha} \Efunc_\mathrm{lw}(w)+ \left( \Efunc_P(u_P) - \Efunc(S_\mathrm{lw}w) \right)\\
&\to& -\mu m(0) + \tfrac{1}{2}\mu^{1+(p-1)\alpha} \Efunc_\mathrm{lw}(w)
\end{eqnarray*} 
as $P \to \infty$, where
$$u_P = \sum_{j \in \Z} (S_\mathrm{lw}w)(\cdot+jP).$$
The result follows by setting $I_\star := -\frac{1}{4}\Efunc_\mathrm{lw}(w)$ and choosing $P_\mu$ large enough so that $\mu^\beta P \geq c_w$
and $| \Efunc_P(u_P) - \Efunc(S_\mathrm{lw}w) | <  \tfrac{1}{4}\mu^{1+(p-1)\alpha} |\Efunc_\mathrm{lw}(w)|$ for $P \geq P_\mu$ (see Lemma~\ref{lemma:PE-convergence}(ii)). \qed

Let us now return to our study of minimisers $\bar{u}_P$ of
$\Efunc_{P,\varrho}$ over $V_{P,\mu}$, which in view of Corollary~\ref{cor:umu} satisfy
\begin{equation}\label{eq:ubar estimate 1}
\Efunc_{P,\varrho}(\bar{u}_P) < -\mu m(0) - \mu^{1+(p-1)\alpha} I_\star
\end{equation}
and of course
\begin{equation} \label{eq:penalisedELeqn}
\mathrm{d}\Efunc_{P,\varrho}[\bar u_{P} ]+\nu_P\mathrm{d}\Qfunc_P[\bar u_{P} ]=0
\end{equation}
for some constant $\nu_P \in \R$, that is
$$
\int_{-\frac{P}{2}}^{\frac{P}{2}}  \left( L\bar u_P + n(\bar u_P)\right)v \dx - 2\varrho^\prime\big( \|\bar u_P\|^2_{H_P^1} \big) \int_{-\frac{P}{2}}^{\frac{P}{2}} \left( \bar u_P v + \bar u_P^\prime v^\prime \right)\dx = \nu_P\int_{-\frac{P}{2}}^{\frac{P}{2}} \bar u_P v \dx
$$
for all $v \in H_P^1$. This equation implies that $\bar u_P^{\prime\prime}$ exists if $\varrho^\prime \big( \|\bar u_P\|^2_{H_P^1} \big)>0$
and that $\bar u_P$ satisfies the equation
\begin{equation}\label{eq:ubar estimate 2}
\nu_P\bar u_P = L\bar u_P + n(\bar u_P) - 2\varrho^\prime \big( \|\bar u_P\|^2_{H_P^1} \big)\left( \bar u_P - \bar u_P^{\prime\prime} \right).
\end{equation}

\begin{lemma}\label{lemma:penmin a}
The estimate
$$\nu_P - m(0) > \tfrac{1}{2}I_\star(p+1)\mu^{(p-1)\alpha}+\bigO(\|\bar{u}_P\|_\infty^{p+\delta-1})-c_\varrho\mu^{1+\varepsilon}$$
holds uniformly over the set of minimisers $\bar{u}_P$ of
$\Efunc_{P,\varrho}$ over $V_{P,\mu}$ and $P \geq P_\mu$.
Here $\varepsilon$ is a positive constant and $c_\varrho$ vanishes when $\varrho=0$.
\end{lemma}
\proof In this proof all estimates hold uniformly in $P \geq P_\mu$.

Inequality~\eref{eq:ubar estimate 1} asserts that
$$
-\int_{-\frac{P}{2}}^{\frac{P}{2}} N(\bar u_P)\dx - \frac{1}{2} \int_{-\frac{P}{2}}^{\frac{P}{2}} \bar u_P L \bar u_P\dx + \varrho\big( \|\bar u_P\|^2_{H_P^1} \big) < - m(0) \mu,
$$
for all $P \geq P_\mu$, and assumption~\eref{eq:mmax} implies that
\begin{equation}\label{eq:m0mu}
\frac{1}{2} \int_{-\frac{P}{2}}^{\frac{P}{2}} \bar u_P L \bar u_P\dx \leq \frac{m(0)}{2} \int_{-\frac{P}{2}}^{\frac{P}{2}} \bar u_P^2\dx = m(0) \mu.
\end{equation}
Adding these inequalities, we find that
$$
\varrho\big( \|\bar u_P\|^2_{H_P^1} \big) \leq \int_{-\frac{P}{2}}^{\frac{P}{2}} N(\bar u_P)\dx
\leq c \|\bar u_P\|_{\infty}^{p-1} \int_{-\frac{P}{2}}^{\frac{P}{2}} \bar u_P^2\dx \leq c \mu^{(p+3)/4},
$$
where we have estimated $\|\bar u_P\|_\infty \leq c \mu^\frac{1}{4}$ (see~\eref{eq:helpfulest}).
Using property~(iii) of the penalisation function, we conclude that
\begin{equation}\label{eq:rho estimates}
\varrho^\prime \big( \|\bar u_P\|^2_{H_P^1} \big) \leq c \mu^{1+\varepsilon}.
\end{equation}

Multiplying~\eref{eq:ubar estimate 2} by $\bar u_P$ and integrating over $(-\frac{P}{2},\frac{P}{2})$,
one finds that
\begin{eqnarray*}
\afl 2 \nu_P\mu & = & (p+1)\int_{-\frac{P}{2}}^{\frac{P}{2}} \left( \frac{1}{2}\, \bar u_P L \bar u_P + N(\bar u_P) \right)\dx - \frac{p-1}{2} \int_{-\frac{P}{2}}^{\frac{P}{2}}   \bar u_P L \bar u_P\dx \\
\afl & & \quad\mbox{} - \int_{\frac{P}{2}}^{\frac{P}{2}} \left( (p+1)N(\bar u_P) - \bar u_P n(\bar u_P) \right) \dx - 2 \varrho^\prime\big( \|\bar u_P \|_{H_P^1}^2 \big) \, \|\bar u_P \|_{H_P^1}^2 \\
\afl & = &  -(p+1)\Efunc_{P,\varrho}(\bar u_P)  - \frac{p-1}{2} 
\int_{-\frac{P}{2}}^{\frac{P}{2}}   \bar u_P L \bar u_P\dx    + (p+1) \varrho \big( \|\bar u_P \|_{H_P^1}^2 \big) \\ 
\afl & &  \quad\mbox{} + \bigO(\|\bar u_P\|_\infty^{p+\delta-1}\|\bar u_P\|_{L_P^2}^2)
- 2 \varrho^\prime\big( \|\bar u_P \|_{H_P^1}^2 \big) \, \|\bar u_P \|_{H_P^1}^2
\end{eqnarray*}
because $(p+1)N(u(x))-un(u(x))=\bigO(|u(x)|^{p+\delta+1})$ uniformly over $u \in V_P$ and $x \in \R$. It follows that
$$
\nu_P > m(0) + \tfrac{1}{2}I_\star(p+1) \mu^{(p-1)\alpha}  + O(\|\bar{u}_P\|_\infty^{p+\delta-1})- c \mu^{1+\varepsilon},
$$
where we have used inequalities~\eref{eq:ubar estimate 1},~\eref{eq:m0mu} and~\eref{eq:rho estimates}.\qed

It follows from Lemma~\ref{lemma:penmin a} and the estimate $\|\bar u_P\|_\infty \leq c \mu^\frac{1}{4}$ (see~\eref{eq:helpfulest})
that 
$\nu_P > \frac{3}{4}m(0)$ uniformly over the set of minimisers $\bar{u}_P$ of
$\Efunc_{P,\varrho}$ over $V_{P,\mu}$ and $P \geq P_\mu$. This bound is used in the following estimate
of the size of $\bar{u}_P$.

\begin{lemma} \label{lemma:penmin b}
The estimate
$$\|\bar{u}_P\|_1^2\le c\mu$$
holds uniformly over the set of minimisers of $\Efunc_{P,\varrho}$ over $V_{P,\mu}$
and $P \geq P_\mu$.
\end{lemma}
\proof In this proof all estimates again hold uniformly in $P \geq P_\mu$.

Multiplying~\eref{eq:ubar estimate 2} by $\bar u_P- \bar u_P^{\prime\prime}$ if $\varrho^\prime \big( \|\bar u_P\|^2_{H_P^1} \big)>0$ or
applying the operator $\bar u_P+\bar u_P^\prime\frac{\mathrm{d}}{\mathrm{d}x}$ if $\varrho^\prime \big( \|\bar u_P\|^2_{H_P^1} \big)=0$, we find that
\begin{eqnarray*}
\afl \nu_P \| \bar u_P \|_{H_P^1}^2 & = & \int_{-\frac{P}{2}}^{\frac{P}{2}} 
\left( \bar u_P L \bar u_P +
\bar u_P^\prime L \bar u_P^\prime\right) \dx 
+ \int_{-\frac{P}{2}}^{\frac{P}{2}} \left( \bar u_P n(\bar u_P) + | \bar u_P^\prime |^2 n^\prime(\bar u_P) \right)\dx\nonumber \\
\afl & & \quad\mbox{} - 2 \varrho^\prime\big( \|\bar u_P \|_{H_P^1}^2 \big) \left( \| \bar u_P \|_{H_P^2}^2 + 2 \int_{-\frac{P}{2}}^{\frac{P}{2}} \left| \bar u_P^\prime \right|^2\dx \right)\nonumber \\
\afl & \leq &  c \| \bar u_P \|_{H_P^{1+\frac{m_0}{2}}}^2 + \left( \sup_{|x| \leq  \|\bar u_P\|_\infty } \!\!\!\!\!|n^\prime(x)|\right)  \|\bar u_P \|_{H_P^1}^2
\end{eqnarray*}
because $m \in S_\infty^{m_0}(\R)$ and
$$|n(u_P(x))| \leq \left(\sup_{|x| \leq  \|\bar u_P\|_\infty } \!\!\!\!\!|n^\prime(x)|\right) |u_P(x)|$$
uniformly over $x \in \R$.
Because $\sup_{|x| \leq  \|\bar u_P\|_\infty } |n^\prime(x)| \to 0$ as $\|\bar u_P\|_\infty \to 0$ and hence as $\|\bar u_P\|_{H_P^1} \to 0$ this quantity is bounded by $\frac{3}{4}m(0)$ for sufficiently small values of $R$, so that
$$
\| \bar u_P \|_{H_P^1}^2 \leq c\|\bar u_P \|_{H_P^{1+\frac{m_0}{2}}}^2.
$$
Estimating
$$
\| \bar u_P \|_{H_P^{1+\frac{m_0}{2}}}^2 \leq
\left\{\begin{array}{ll}
\| \bar u_P \|_{L_P^2}^{|m_0|} \| \, \bar u_P \|_{H_P^1}^{2-|m_0|}, & |m_0| < 2, \\[2mm]
\| \bar u_P \|_{L_P^2}^2, & |m_0| \geq 2 \end{array}\right.
$$
shows that
$$
\| \bar u_P \|_{H_P^1}^2 \leq c\| \bar u_P \|_{L_P^2}^2 \leq c\mu.\eqno{\Box}
$$


\begin{theorem}[Existence of periodic minimisers]\label{theorem:P-minimiser}
For each $P \geq P_\mu$ there exists a function $\bar u_P \in U_{P,\mu}$ which minimises $\Efunc_P$ over $U_{P,\mu}$, so that
$$\Efunc_P(\bar u_P) = I_{P,\mu}:= \inf \left\{ \Efunc_{P}(u) \colon u \in U_{P,\mu} \right\},$$
and satisfies the Euler-Lagrange equation
$$\Efunc_P^\prime(\bar u_P) + \nu_P \Qfunc_P^\prime(\bar u_P)=0$$
for some real number $\nu_P$; it is therefore
a periodic solution of the travelling-wave equation~\eref{eq:steadyproblem} with wave speed $\nu_P$.
Furthermore
$$\| \bar u_P \|_{H_P^1}^2 \leq c\mu, \qquad 0 < \nu_P \leq c$$
uniformly over $P \geq P_\mu$.

\end{theorem}
\proof Let $\bar{u}_P$ be a minimiser of $\Efunc_{P,\varrho}$ over $V_{P,\mu}$. It follows from
Lemma~\ref{lemma:penmin b} that $\|\bar{u}_P\|_1^2 \leq c\mu$, so that $\varrho(\bar u_P)$ and
$\varrho^\prime(\bar u_P)$ vanish.  In particular, $\bar u_P$ belongs to $U_{P,\mu}$, and since it minimises $\Efunc_{P,\varrho}$ over $V_{P,\mu}$ it certainly
minimises $\Efunc_{P,\varrho}=\Efunc_P$ over $U_{P,\mu}$. Furthermore, equation~\eref{eq:penalisedELeqn}
is equivalent to
$$\Efunc_P^\prime(\bar u_P) + \nu_P\Qfunc_P^\prime (\bar u_P)=0,$$
from which it follows that
$$\nu_P = -\frac{1}{2\mu}\langle \Efunc_P^\prime(\bar u_P),\bar u_P \rangle_{L_P^2} \leq \frac{c}{\mu}\|\bar{u}_P\|_{H_P^1}^2 \leq c.\eqno{\Box}$$


\subsection*{Construction of a special minimising sequence for $\Efunc$}
We proceed by extending the minimisers $\bar u_P$ of $\Efunc_P$ over $U_{P,\mu}$ found above
to functions in $H^1(\R)$ by scaling, translation and truncation in the following manner. For each sufficiently large value of $P$ there exists an open subinterval\linebreak
$I_P := (x_P - \frac{1}{2} P^\frac{1}{4}, x_P + \frac{1}{2} P^\frac{1}{4})$ of
$(-\frac{P}{8}, \frac{P}{8})$ such that $\| \bar u_P \|_{H^1(I_P)} < P^{-\frac{1}{4}}$; we may assume that this property holds for all $P \geq P_\mu$.
Let $\chi \colon [0,\infty) \to [0,\infty)$ be a smooth, increasing `cut-off' function with
\begin{equation*}\label{eq:chi}
\chi(r) := 
\left\{\begin{array}{ll}
0, \qquad & 0 \leq r \leq 1/2,\\[1mm]
1, \qquad & r \geq 1,
\end{array}\right.,
\end{equation*}
let $u_P$ be the $P$-periodic function defined by
$$
u_P(x) :=  A_P v_P \left(x+\tfrac{P}{2} \right),
$$
where
$$
v_P(x)|_{[-\frac{P}{2},\frac{P}{2}]} := \chi\left( \frac{2|x|}{P^\frac{1}{4}}\right) \bar u_P(x + x_P),\qquad A_P = \frac{\sqrt{2\mu}}{\|v_P\|_{L_P^2}},
$$
and finally define $\tilde u_P \in H^1(\R)$ by the formula
$$
\tilde u_P(x) :=
\left\{\begin{array}{ll}
u_P(x), \qquad &|x| \leq \frac{P}{2},\\[1mm]
0, 	\qquad &|x| > \frac{P}{2},
\end{array}\right.
$$
so that
$$u_P = \sum_{j \in \Z} \tilde{u}_P(\cdot+jP).$$

Let us examine the sequence $\{\tilde u_n\}_{n \in \N_0}$, where $\tilde u_n := \tilde u_{P_n}$ and $\{P_n\}_{n\in\N_0}$ is an increasing, unbounded sequence of positive real numbers with $P_0 \geq P_\mu$.

\begin{theorem}[Special minimising sequence for $\Efunc$]\label{theorem:minimisers}
The sequence $\{\tilde u_n\}_{n \in \N_0}$ is a minimising sequence for $\Efunc$ over $U_{\mu}$
which satisfies
$$
\sup_{n \in \N_0} \| \tilde u_n \|_1^2 \leq c \mu, \qquad
\lim_{n\to\infty}\|\Efunc^\prime(\tilde u_n)+\nu_n \Qfunc^\prime(\tilde u_n)\|_1=0,
$$
where $\nu_n=\nu_{P_n}$, $n \in \N_0$.
\end{theorem}
\proof
Observe that
\begin{eqnarray*}
\lefteqn{\left\| u_P - \bar u_P\left( \cdot + x_P +  \tfrac{P}{2}\right) \right\|_{L_P^2}^2}\qquad\\ 
&= &\int_{-\frac{P}{2}}^{\frac{P}{2}} \left| A_P \, \chi\left( \frac{2|x|}{P^\frac{1}{4}} \right) - 1 \right|^2 |\bar u_P(x + x_P)|^2\dx\\
&= &\int_{|x| < \frac{1}{2} P^\frac{1}{4}}  \left| A_P \, \chi\left( \frac{2|x|}{P^\frac{1}{4}} \right) - 1 \right|^2 |\bar u_P(x + x_P)|^2\dx\\ 
&& \quad\mbox{} + |A_P - 1|^2 \underbrace{\int_{|x| > \frac{1}{2} P^\frac{1}{4}} |\bar u_P(x + x_P)|^2\dx}_{\displaystyle \leq \|\bar u_P\|_{H_P^1}^2 < R} \\
& \to & 0
\end{eqnarray*}
as $P \to \infty$; the first integral vanishes by the choice of the intervals $I_P$, while the factor $A_P - 1$ also vanishes
because $\lim_{P\to \infty} \| v_P \|_{L_P^2} = \lim_{P \to \infty} \|\bar u_P \|_{L_P^2} = \sqrt{2\mu}$. Similarly,
\begin{eqnarray*}
\lefteqn{\left\| u_P^\prime - \bar u_P^\prime\left( \cdot + x_P +  \tfrac{P}{2}\right) \right\|_{L_P^2}^2}\qquad\\ 
&= &\int_{-\frac{P}{2}}^{\frac{P}{2}} \left| A_P \, \chi\left( \frac{2|x|}{P^\frac{1}{4}} \right) - 1 \right|^2 |\bar u_P^\prime(x + x_P)|^2\dx\\
& & \quad\mbox{}+\frac{4}{P^\frac{1}{2}}\int_{-\frac{P}{2}}^{\frac{P}{2}} \left| A_P \, \chi^\prime\left( \frac{2|x|}{P^\frac{1}{4}} \right) \right|^2 |\bar u_P(x + x_P)|^2\dx\\
& \to & 0
\end{eqnarray*}
as $P \to \infty$ (the above argument shows that the first integral vanishes, while the second integral is bounded). It follows that
$$
\left\| u_P - \bar u_P\left( \cdot + x_P +\tfrac{P}{2}\right) \right\|_{H_P^1} \to 0 \quad\mbox{ as }P \to \infty,
$$
and this result shows in particular that
$$\|\tilde u_P\|_1 = \|u_P\|_{H_P^1} \leq \left\| u_P - \bar u_P\left( \cdot + x_P +\tfrac{P}{2}\right) \right\|_{H_P^1}  + \| \bar u_P\|_{H_P^1} \leq c \mu$$
for $P \geq P_\mu$ (where $P_\mu$ is replaced with a larger constant if necessary).

Next note that
\begin{eqnarray*}
\Efunc_P(u_P)-\Efunc_P(\bar u_P) & = & \Efunc_P(u_P)-\Efunc_P\left(\bar u_P\left(\cdot+x_P+\tfrac{P}{2}\right)\right) \\
& \leq & \sup_{u \in U_P} \|\Efunc_P^\prime(u)\|_{L_P^2} \left\| u_P - \bar u_P\left( \cdot + x_P +\tfrac{P}{2}\right) \right\|_{L_P^2} \\
& \to & 0
\end{eqnarray*}
as $P \rightarrow \infty$ (because $\|\Efunc_P^\prime(u)\|_{L_P^2}$ is
bounded uniformly over $u \in U_P$ and $P>0$) and 
$$\Efunc(\tilde u_P) - \Efunc_P(u_P) \to 0$$
as $P \to \infty$ (Lemma~\ref{lemma:PE-convergence}(ii)). Observe further that $I_{P,\mu} \rightarrow I_\mu$ as $P \rightarrow \infty$:
\begin{itemize}
\item
Take $\tilde w \in C_\mathrm{c}^\infty(\R)$ with $\Qfunc(\tilde w) = \mu$, so that $w_P := \sum_{j \in \Z} \tilde w( \cdot + jP)$ satisfies $I_{P,\mu} \leq \Efunc_P(w_P)$ and $\Efunc_P(w_P) \to \Efunc(\tilde w)$ as $P \to \infty$ (see Lemma~\ref{lemma:PE-convergence}(ii)). It follows that $\limsup_{P \to \infty} I_{P,\mu} \leq \Efunc(\tilde w)$, and hence that
$$
\limsup_{P \to \infty} I_{P,\mu} \leq \inf\left\{ \Efunc(\tilde u) \colon \tilde u \in C_\mathrm{c}^\infty(\R) \cap U_{\mu} \right\} = I_\mu.
$$
\item
On the other hand, 
\begin{eqnarray*}
I_\mu &\leq &\Efunc(\tilde u_{P})\\
&= &\left( \Efunc(\tilde u_{P}) - \Efunc_{P}(u_{P}) \right)+ \left(\Efunc_{P}(u_{P}) - \Efunc_{P}(\bar u_{P}) \right) + I_{P,\mu},
\end{eqnarray*}
in which the first and second terms on the right-hand side vanish as $P \to \infty$, so that
$$
I_\mu \leq \liminf_{P \to \infty} I_{P,\mu}.
$$
\end{itemize}
We conclude that
$$\Efunc(\tilde{u}_P) = \left(\Efunc(\tilde u_{P}) - \Efunc_{P}(u_{P}) \right)+ \left(\Efunc_{P}(u_{P}) - \Efunc_{P}(\bar u_{P}) \right) + I_{P,\mu} \to I_\mu$$
as $P \to \infty$.

Similarly, note that
\begin{eqnarray*}
\afl\|\Efunc_P^\prime(u_P)-\Efunc_P^\prime(\bar u_P)\|_{H_P^1} & = & \left\| \Efunc_P^\prime(u_P)-\Efunc_P^\prime\left(\bar u_P\left(\cdot+x_P+\tfrac{P}{2}\right)\right)\right\|_{H_P^1} \\
\afl & \leq & \sup_{u \in U_P} \|\mathrm{d}\Efunc_P^\prime[u]\|_{H_P^1 \rightarrow H_P^1} \left\| u_P - \bar u_P\left( \cdot + x_P +\tfrac{P}{2}\right) \right\|_{H_P^1} \\
\afl & \to & 0
\end{eqnarray*}
as $P \rightarrow \infty$ (it follows from
the calculation $\mathrm{d}\Efunc_P^\prime[u](v) = -Lv - n^\prime(u)v$ that
\begin{equation}\label{eq:nmustbeC2}
\afl\|\mathrm{d}\Efunc_P^\prime[u]\|_{H_P^1 \rightarrow H_P^1} \leq c\left(m(0)
+\!\!\!\sup_{|x| \leq c_sR}  |n^\prime(x)| +\!\!\! \sup_{|x| \leq c_sR} |n^{\prime\prime}(x)|\right) \leq c
\end{equation}
uniformly over $u \in U_P$ and $P>0$),
and Lemma~\ref{lemma:PE-convergence}(ii) shows that
$$
\lim_{P \to \infty} \| \Efunc^\prime (\tilde u_P) - \Efunc^\prime_P (u_P) \|_{H^1(-\frac{P}{2},\frac{P}{2})} = 0, \qquad
\lim_{P \to \infty} \|\Efunc^\prime (\tilde u_P) \|_{H^1(\{|x| > \frac{P}{2}\})}=0;
$$
the same results hold for $\Qfunc$, $\Qfunc_P$.
We conclude that
\begin{eqnarray*}
\afl \|\Efunc^\prime(\tilde{u}_P)+\nu_P \Qfunc^\prime(\tilde{u}_P)\|_1 \\
\afl\quad\leq \|\Efunc^\prime(\tilde{u}_P)- \Efunc_P^\prime(u_P)\|_{H^1(-\frac{P}{2},\frac{P}{2})}
+\nu_P  \|\Qfunc^\prime(\tilde{u}_P)- \Qfunc_P^\prime(u_P)\|_{H^1(-\frac{P}{2},\frac{P}{2})}\\
\afl \hspace{1cm}\mbox{}+ \|\Efunc_P^\prime(u_P) - \Efunc_P^\prime(\bar{u}_P)\|_{H^1(-\frac{P}{2},\frac{P}{2})}
+ \nu_P\|\Qfunc_P^\prime(u_P) - \Qfunc_P^\prime(\bar{u}_P)\|_{H^1(-\frac{P}{2},\frac{P}{2})}\\
\afl \hspace{1cm}\mbox{}+ \|\Efunc_P^\prime(\bar{u}_P)+\nu_P \Qfunc_P^\prime(\bar{u}_P)\|_{H^1(-\frac{P}{2},\frac{P}{2})}
+ \|\Efunc^\prime(\tilde{u}_P)\|_{H^1(\{|x|>\frac{P}{2}\})}\\
\afl \hspace{1cm}\mbox{}+\nu_P\|\Qfunc^\prime(\tilde{u}_P)\|_{H^1(\{|x|>\frac{P}{2}\})} \\
\afl\quad \to 0
\end{eqnarray*}
as $P \to \infty$ because $\{\nu_P\}$ is bounded.\qed

\section{Strict subadditivity} \label{sec:SSA}

In this section we show that the quantity
$$I_\mu :=\inf \left\{ \Efunc(u) \colon u \in U_{\mu} \right\}$$
is \emph{strictly subadditive}, that is
$$
I_{\mu_1 + \mu_2} < I_{\mu_1} + I_{\mu_2} \quad \mbox{ whenever }0 < \mu_1,\mu_2  < \mu_1 + \mu_2 < \mu_\star.
$$
This result is needed in Section~\ref{sec:concentration}  below to exclude `dichotomy' when applying the concentration-compactness
principle to a minimising sequence $\{u_n\}_{n \in \N_0}$ for $\Efunc$ over $U_\mu$. It is proved by approximating the nonlinear term $\Nfunc(u_n)$ by its leading-order
homogeneous part $-\int_\R N_{p+1}(u_n)\dx$ (strict subadditivity for a problem with a homogeneous nonlinearity follows by a straightforward scaling argument). However, the requisite estimate
$$
\int_\R  N_\mathrm{r}(u_n)\dx  =o(\mu^{1+(p-1)\alpha})
$$
may not hold
for a general minimising sequence; it does however hold for the special minimising sequence
$\{\tilde{u}_n\}_{n \in \N_0}$ constructed in Section~\ref{sec:periodic} above.




\subsection*{Scaling}

We now examine functions $u \in U_\mu$ which are `near minimisers' of $\Efunc$ in the sense that
\begin{equation} \label{eq:General ms estimates a} 
\Efunc(u)< -\mu m(0)-I_\star\mu^{1+(p-1)\alpha},\qquad
\|\Efunc^\prime(u)+\nu \Qfunc^\prime(u)\|_1 \le c\mu^{N}
\end{equation}
for some $\nu \in \R$ and natural number $N\ge \max\{\frac{1}{2}(1+4j_\star\beta),1+(p-1)\alpha\}$. We show that their low-wavenumber part is a long wave which `scales' in a fashion similar to the
\emph{Ansatz}~\eref{eq:longwaveansatz};
this result allows us to conclude in particular that $\|u\|_\infty \leq c \mu^{\alpha-\varepsilon}$ for any $\varepsilon>0$
(see Corollary~\ref{cor:usupnorm}).

Our results are obtained by studying the identity
\begin{equation}\label{eq:ELeqn}
\nu u - Lu  = n(u) + \Efunc^\prime(u)+\nu \Qfunc^\prime(u);
\end{equation}
they apply to minimisers $u$ of $\Efunc$ over $U_\mu$, for which $\Efunc^\prime(u)+\nu \Qfunc^\prime(u)=0$
for some Lagrange multiplier $\nu$, and to the functions $\tilde{u}_n$ in the minimising sequence $\{\tilde{u}_n\}_{n \in \N_0}$, which
satisfy $\lim_{n \rightarrow \infty} \|\Efunc(\tilde{u}_n) + \nu_n \Qfunc(\tilde{u}_n)\|_1 =0$.
(Without loss of generality we may assume that
$\nu_n$ does not depend upon $n$: the bounded sequence
$\{\nu_n\}_{n \in \N_0}$ has a convergent subsequence whose limit $\nu$ satisfies
$$\lim_{n \rightarrow \infty} \|\Efunc^\prime(\tilde{u}_n)+\nu \Qfunc^\prime(\tilde{u}_n)\|_1 =0$$
because $\{\|\Qfunc^\prime(\tilde{u}_n)\|_1\}_{n \in \N_0}$ is bounded.)

We begin with the following preliminary result, which is proved in the same fashion as
Lemma~\ref{lemma:penmin a}.

\begin{proposition} \label{prop:speed}
The estimate
$$\nu - m(0) > \tfrac{1}{2}I_\star(p+1)\mu^{(p-1)\alpha}+\bigO(\|u\|_\infty^{p+\delta-1})+\bigO(\mu^{N-\frac{1}{2}})$$
holds uniformly over the set of $u \in U_\mu$ satisfying~\eref{eq:General ms estimates a}.
\end{proposition}

According to Proposition~\ref{prop:speed} one may replace~\eref{eq:General ms estimates a} by
\begin{equation}
\nu - m(0) > \bigO(\|u\|_\infty^{p+\delta-1}), \qquad \|\Efunc^\prime(u)+\nu \Qfunc^\prime(u)\|_1 \le c\mu^{N}, \label{eq:General ms estimates b}
\end{equation}
and most of the results in the present section apply to this more general situation. In particular,
estimating
$$\|u\|_\infty \leq c \|u\|_0^\frac{1}{2}\|u\|_1^\frac{1}{2} \leq c\mu^\frac{1}{4},\qquad u \in U_\mu$$
we find that $\nu>\frac{3}{4}m(0)$; our next result is obtained from this bound in the same fashion as Lemma~\ref{lemma:penmin b}.

\begin{proposition}
The estimate
$$\|u\|_1^2\le c\mu$$
holds uniformly over the set of $u \in U_\mu$ satisfying~\eref{eq:General ms estimates b}.
\end{proposition}

The next step is to decompose a function $u \in H^1(\R)$
into low- and high-wavenumber parts in the following manner. Choose $k_0>0$ so that $m(k) \leq \frac{1}{2}m(0)$
for $|k| \geq k_0$, let $\xi$ be the characteristic function of the set $[-k_0,k_0]$, and write $u=u_1+u_2$, where 
$$
\hat u_1(k) := \xi(k)\hat u(k), \qquad \hat u_2(k) := (1-\xi(k))\hat u(k). 
$$
We proceed by writing~\eref{eq:ELeqn}
as coupled equations for the low- and high-wavenumber parts of $u$, namely
\begin{eqnarray}
(\nu-m)\hat{u}_1 = \xi \FF[ n(u)+\Efunc^\prime(u)+\nu\Qfunc^\prime(u)], \label{eq:Euler-Lagrange 1 gen} \\
(\nu-m) \hat{u}_2 = (1-\xi) \FF[ n(u)+\Efunc^\prime(u)+\nu\Qfunc^\prime(u)], \label{eq:Euler-Lagrange 2 gen}
\end{eqnarray}
and estimating $u_1$ using the weighted norm
\begin{equation}\label{eq:weightednorm}
\nn v\nn_{\tau,\mu} := \left( \int_\R \left( v^2 + \mu^{-4j_\star\tau \beta} (v^{(2j_\star)})^2 \right)\dx \right)^{\!\!\frac{1}{2}},\qquad \tau<1
\end{equation}
for $H^{2j_\star}(\R)$, which is useful in estimating the $L^\infty(\R)$-norm of $u_1$ and its derivatives.

\begin{proposition}\label{prop:supnorm}
The estimate
$$\|v^{(j)}\|_\infty \leq c\mu^{(j+\frac{1}{2})\tau\beta}\nn v \nn_{\tau,\mu}, \qquad j=0,\ldots,2j_\star-1$$
holds for all $v \in H^{2j_\star}(\R)$.
\end{proposition}
{\bf Proof.} Observe that
$$
\|v^{(j)}\|_\infty^2\! \leq\! \frac{1}{2\pi}\|k^{j}\hat v\|_{L^1(\R)}^2 \leq \frac{1}{2\pi}\!\!\left(\int_\R\!\frac{k^{2j}}{1+\mu^{-4j_\star\tau \beta}k^{4j_\star}} \dk\!\!\right)\! \nn v\nn_{\tau,\mu}^2
\leq c\mu^{(2j+1)\tau \beta}\nn v\nn_{\tau,\mu}^2. \eqno{\Box}
$$
\medskip


\begin{theorem}[Scaling]\label{theorem:scaling gen}
Choose $\tau <1$. The estimates
$$
\nn u_1\nn_{\tau,\mu}^2 \leq  c_\tau\mu \qquad \|u_{2}\|_1^2 \leq c_\tau \mu^{\tau\beta(p-1)+p}
$$
hold for all $u \in U_\mu$ which satisfy~\eref{eq:General ms estimates b}.
\end{theorem}
\proof Note that
$\nu-m(k) \geq \frac{1}{4}m(0)$ for $|k| \geq k_0$ (since $\nu>\frac{3}{4}m(0)$), so that
$$
\FF^{-1}[(\nu-m)^{-1} (1-\xi)\FF(\cdot)] \in B(H^1(\R),H^1(\R)),
$$
where the operator norm is bounded uniformly over $\nu>\frac{3}{4}m(0)$,
and it follows from equation~\eref{eq:Euler-Lagrange 2 gen} that
\begin{eqnarray*}
\|u_2\|_1 & \leq & c (\|n(u)\|_1+\|\Efunc^\prime(u)+\nu \Qfunc^\prime(u)\|_1)\\  
&\le&
c( \mu^\frac{1}{2} \|u_1\|_\infty^{p-1}+\mu^{\frac{1}{2}(p-1)}\|u_2\|_1+ \mu^N),
\end{eqnarray*}
where we have estimated
\begin{eqnarray*}
\|n(u)\|_1^2 & = & \|n(u)\|_0^2 + \|n^\prime(u)u^\prime\|_0^2 \\
& \leq & c \|u\|_\infty^{2p-2} \|u\|_1^2 \\
& \leq & c (\|u_1\|_\infty^{2p-2} \|u_1\|_1^2 + \|u_2\|_\infty^{2p-2} \|u_1\|_1^2 +\|u\|_\infty^{2p-2} \|u_2\|_1^2) \\
& \leq & c (\|u_1\|_\infty^{2p-2} \|u\|_1^2 + \|u_2\|_1^{2p-2} \|u_1\|_1^2 +\|u\|_1^{2p-2} \|u_2\|_1^2) \\
& \leq & c (\|u_1\|_\infty^{2p-2} \|u\|_1^2 +\|u\|_1^{2p-2} \|u_2\|_1^2) \\
& \leq & c (\mu \|u_1\|_\infty^{2p-2} + \mu^{p-1}\|u_2\|_1^2).
\end{eqnarray*}
We conclude that
\begin{equation}\label{eq:u2mu gen}
\|u_2\|_1 \le c(\mu^\frac{1}{2}  \|u_1\|_{\infty}^{p-1}+\mu^N).
\end{equation}

Turning to equation~\eref{eq:Euler-Lagrange 1 gen}, observe that
$$\nu - m(k)  > (\nu-m(0)) - \frac{cm^{(2j_\star)}(0)}{(2j_\star)!} k^{2j_\star} > - \frac{cm^{(2j_\star)}(0)}{(2j_\star)!} k^{2j_\star}+\bigO(\|u\|_\infty^{p+\delta-1})$$ for $|k|<k_0$ and uniformly over $u \in U_\mu$, so that
\begin{eqnarray*}
\afl \int_\R |u_1^{(2j_\star)}|^2\dx \\
\afl\qquad\leq c \int_\R ( \nu - m(k) )^2 |\hat u_1(k)|^2\dk+c\|u_1\|_1^2\|u\|_\infty^{2(p-1)}\\ 
\afl\qquad\leq c (\|n(u)\|_0^2+\|\Efunc^\prime(u)+\nu \Qfunc^\prime(u)\|_0^2
+\|u_1\|_1^2\|u_1\|_\infty^{2(p-1)}+\|u_1\|_1^2\|u_2\|_1^{2(p-1)}) \\
\afl\qquad\leq c (\|n(u)\|_0^2+\|\Efunc^\prime(u)+\nu \Qfunc^\prime(u)\|_0^2
+\|u\|_1^2\|u_1\|_\infty^{2(p-1)}+\|u\|_1^{2(p-1)}\|u_2\|_1^2) \\
\afl\qquad\leq c (\|n(u)\|_0^2+\|\Efunc^\prime(u)+\nu \Qfunc^\prime(u)\|_0^2
+\mu\|u_1\|_\infty^{2(p-1)}+\|u_2\|_1^2) \\
\afl\qquad\leq c(\mu \| u_1 \|_\infty^{2(p-1)} +\|u_2\|_1^2+\mu^{2N})\\
\afl\qquad\leq c(\mu \| u_1 \|_\infty^{2(p-1)}+\mu^{2N}) \\
\afl\qquad\leq c (\mu^{1+(p-1)\tau\beta}\nn u_1\nn_{\tau,\mu}^{2(p-1)} +\mu^{2N})\\
\afl\qquad\leq c\left( \mu^{1+(p-1)(\tau\beta+1)}\left(\frac{\nn u_1\nn_{\tau,\mu}}{\mu^\frac{1}{2}}\right)^{\!\!2(p-1)}+\mu^{2N}\right),
\end{eqnarray*}
where we have estimated
$$\|n(u)\|_0^2 \leq c (\|u_1\|_\infty^{2p-2} \|u_1\|_0^2 + \|u_2\|_1^{2p}) \leq c (\mu\|u_1\|_\infty^{2p-2}  + \|u_2\|_1^2)$$
and used~\eref{eq:u2mu gen} and Proposition~\ref{prop:supnorm}. Multiplying this estimate by $\mu^{-4j_\star\tau\beta}$
and adding the inequality
$\int_\R u_1^2\dx \leq 2\mu$, one finds that
$$
\nn u_1 \nn_{\tau,\mu}^2 \leq c  \mu \left(1+ \mu^{(1-\tau)(p-1)}\left(\frac{\nn u_1\nn_{\tau,\mu}^2}{\mu}\right)^{\!\!p-1}\right).
$$

Define $Q=\{\tau \in (-\infty,1): \nn u_1 \nn_{\tau,\mu}^2 \leq c_\tau \mu\}$. The inequality $\nn u_1 \nn_{\tau_1,\mu}^2
\leq \nn u_1 \nn_{\tau_2,\mu}^2$ for $\tau_1 \leq \tau_2$ shows that $(-\infty,\tau] \subset Q$
whenever $\tau \in Q$; furthermore $(-\infty,0] \subseteq Q$ because $\nn u_1 \nn_{0,\mu}^2 \leq \| u_1 \|_0^2 \leq 2 \mu$.
Suppose that $\tau_\star:=\sup Q$ is strictly less than unity, choose $\varepsilon>0$ so that
$\tau_\star+(1+8j_\star\beta)\varepsilon<1$ and observe that
\[
\frac{\nn u_1 \nn_{\tau_\star+\varepsilon,\mu}^2}{\mu}\ \leq\ c\Bigg(1+\mu^{(1-\tau_\star-(1+8j_\star\beta)\varepsilon)(p-1)}
\underbrace{\left(\frac{\nn u_1 \nn_{\tau_\star-\varepsilon,\mu}^2}{\mu}\right)^{\!\!p-1}}_{\displaystyle \leq c_{\tau^\star-\varepsilon}}\Bigg),
\]
which leads to the contradiction that $\tau_\star+\varepsilon \in Q$. It follows that $\tau_\star=1$ and
$\nn u_1 \nn_{\tau,\mu}^2 \leq c_\tau\mu$ for each $\tau<1$. 

The bound for $\|u_2\|_1^2$ follows from inequality~\eref{eq:u2mu gen}, Proposition~\ref{prop:supnorm} and the bound for
$\nn u_1 \nn_{\tau,\mu}^2$.\qed

\begin{corollary}\label{cor:usupnorm}
Choose $\tau<1$. The estimate
$$\|u\|_\infty \leq c_\tau \mu^\alpha\mu^{(1-\tau)(\frac{1}{2}-\alpha)}$$
holds for all $u \in U_\mu$ which satisfy~\eref{eq:General ms estimates b}.
\end{corollary}
\proof
Using Proposition~\ref{prop:supnorm}, Theorem~\ref{theorem:scaling gen} and the relation $\beta=2\alpha-1$, one finds that
\begin{eqnarray*}
\|u\|_\infty & \leq & \|u_1\|_\infty + \|u_2\|_\infty \\
& \leq & c (\mu^{\frac{\tau\beta}{2}}\nn u_1 \nn_{\tau,\mu} + \|u_2\|_1) \\
& \leq & c_\tau(\mu^{\frac{1}{2}+\frac{\tau\beta}{2}} + \mu^{\frac{p}{2}+\frac{\tau\beta}{2}(p-1)})\\
& \leq & c_\tau \mu^{\frac{1}{2}+\frac{\tau\beta}{2}} \\
& = & c_\tau \mu^\alpha\mu^{(1-\tau)(\frac{1}{2}-\alpha)}.\\[-6mm]
\end{eqnarray*}
\qed

\begin{corollary} \label{cor:supercritical}
Any function $u \in U_\mu$ satisfying~\eref{eq:General ms estimates a} has the property that
$$\nu-m(0)>0.$$
\end{corollary}
\proof Using Corollary~\ref{cor:usupnorm}, we find that
\begin{eqnarray*}
\|u\|_\infty ^{p+\delta-1} & \leq & c\|u\|_1^{\delta} \|u\|_\infty^{p-1} \\
& \leq & c_\tau  \mu^{\frac{\delta}{2}}\mu^{(p-1)\alpha}\mu^{(1-\tau)(\frac{1}{2}-\alpha)(p-1)}\\
& = & c_\tau \mu^{\frac{\delta}{2}+(1-\tau)(\frac{1}{2}-\alpha)(p-1)} \mu^{(p-1)\alpha} \\
& = & o(\mu^{(p-1)\alpha})
\end{eqnarray*}
uniformly over $u \in D_\mu$ for $\tau$ sufficiently close to $1$, whereby Proposition~\ref{prop:speed} shows that
$$\nu - m(0) > \tfrac{1}{2}I_\star(p+1)\mu^{(p-1)\alpha}+o(\mu^{(p-1)\alpha}) > 0.\eqno{\Box}$$

\subsection*{Strict subhomogeneity}

A function $\mu \mapsto I_\mu$ is said to be \emph{strictly subhomogeneous} on an interval $(0,\mu_\star)$ if
$$
I_{a\mu} < a I_\mu \quad\mbox{ whenever }  \quad 0 < \mu < a\mu  < \mu_\star;
$$  
a straightforward argument shows that strict subhomogeneity implies strict subadditivity on the same interval
(see Buffoni \cite[p.\ 48]{Buffoni04a}).

\begin{proposition}\hspace{1in}\label{prop:nneg}
\begin{itemize}
\item[(i)] Any function $u \in U_\mu$ with the property
\begin{equation}
\Efunc(u)< -\mu m(0)-I_\star\mu^{1+(p-1)\alpha} \label{eq:minest}
\end{equation}
satisfies
$$
\Nfunc(u) \leq -c \mu^{1+(p-1)\alpha}.
$$
This result holds in particular for any minimising sequence 
$\{u_n\}_{n \in \N_0}$ for $\Efunc$ over $U_\mu$.
\item[(ii)]
Any function $u \in U_\mu$ with the property~\eref{eq:minest} satisfies
$$\int_\R  N_{p+1}(u)\dx \geq c \mu^{1+(p-1)\alpha}.$$
This result holds in particular for the sequence $\{\tilde u_n\}_{n \in \N_0}$.
\end{itemize}
\end{proposition}
\proof The first result is a consequence of the equation
$$\Nfunc(u) =  \Efunc(u) - \Lfunc(u)$$
and the estimates~\eref{eq:minest} and
$$-\Lfunc(u) = \frac{1}{2}\int_ \R u Lu \dx \leq \frac{m(0)}{2} \int_\R u^2 \dx = \mu m(0),$$
while the second is obtained from the first using the estimate
\begin{eqnarray}
\left|\int_\R N_\mathrm{r}(u)\dx\right|  & \leq & c\|u\|_1^{2+\delta} \|u\|_\infty^{p-1} \nonumber \\
& \leq & c_\tau \mu^{1+\frac{\delta}{2}}\mu^{(p-1)\alpha}\mu^{(1-\tau)(\frac{1}{2}-\alpha)(p-1)} \nonumber \\
& = & c_\tau \mu^{\frac{\delta}{2}+(1-\tau)(\frac{1}{2}-\alpha)(p-1)} \mu^{1+(p-1)\alpha} \nonumber \\
& = & o(\mu^{1+(p-1)\alpha}) \label{Estimate for Nr}
\end{eqnarray}
for $\tau$ sufficiently close to $1$ (see Corollary~\ref{cor:usupnorm}).\qed


\begin{lemma}\label{lemma:subadditive}
The map $\mu \mapsto I_\mu$  is strictly subhomogeneous for $\mu \in (0,\mu_\star)$.
\end{lemma}
\proof Fix $a > 1$ and note that $\| a^\frac{1}{2} \tilde u_n \|_1^2 \leq c a \mu < R$. We have that
\begin{eqnarray}
I_{a\mu} &\leq& \Efunc( a^\frac{1}{2} \tilde u_n ) \nonumber \\
&=& \Lfunc( a^\frac{1}{2} \tilde u_n ) - \int_\R N_{p+1}( a^\frac{1}{2} \tilde u_n)\dx -
\int_\R N_\mathrm{r}(a \tilde u_n)\dx \nonumber \\
&=& a \Lfunc( \tilde u_n ) - a^{\frac{1}{2}(p+1)} \int_\R N_{p+1}(\tilde u_n)\dx + o(\mu^{1+(p-1)\alpha}) \nonumber \\
&=& a \Efunc( \tilde u_n ) - (  a^{\frac{1}{2}(p+1)} -a) \int_\R N_{p+1}(\tilde u_n)\dx + o(\mu^{1+(p-1)\alpha}) \nonumber \\
& \leq & a \Efunc( \tilde u_n ) - c(  a^{\frac{1}{2}(p+1)} -a) \mu^{1+(p-1)\alpha} +o(\mu^{1+(p-1)\alpha}), \label{eq:ssh2}
\end{eqnarray}
in which we have used Proposition~\ref{prop:nneg}(ii) and the estimate
$$
\left|\int_\R N_\mathrm{r}(a\tilde u_n)\dx\right| \leq ca^{p+\delta+1}\|\tilde u_n\|_1^{2+\delta} \|\tilde u_n\|_\infty^{p-1}
=o(\mu^{1+(p-1)\alpha})
$$
(cf.\ calculation~\eref{Estimate for Nr}).
In the limit $n \rightarrow \infty$ inequality~\eref{eq:ssh2} yields
$$
I_{a\mu} \leq a I_\mu - c(  a^{\frac{1}{2}(p+1)} -a) \mu^{1+(p-1)\alpha} +o(\mu^{1+(p-1)\alpha}),
$$
from which it follows that $I_{a\mu} < a I_\mu$.\qed




\section{Concentration-compactness}\label{sec:concentration} 

In this section we present the proof of Theorem~\ref{theorem:main} with the help of the concentration-compactness principle
(Lions \cite{Lions84a}), which we now recall in a form suitable for our purposes.


\begin{theorem}[Concentration-compactness]\label{theorem:concentration}
Any sequence $\{e_n\}_{n \in \N_0} \subset L^1(\R)$ of non-negative functions with the property that
$$
\int_\R e_n\dx = l > 0
$$
admits a subsequence, denoted again  by $\{e_n\}_{n \in \N_0}$, for which one of the following phenomena occurs.\\[-6pt]

\noindent \underline{Vanishing:}\, For each $r > 0$ one has that 
\begin{equation}\label{eq:vanishing}
\lim_{n\to \infty} \left( \sup_{x_0 \in \R} \int_{B_r(x_0)} e_n\dx \right) = 0.\\[5pt]
\end{equation}

\noindent \underline{Concentration:}\, There is a sequence $\{x_n\}_{n\in \N_0} \subset \R$ with the property that for each $\varepsilon >0$ there exists $r > 0$ with
\begin{equation}\label{eq:concentration}
\int_{B_r(x_n)} e_n\dx \geq l - \varepsilon,
\end{equation}
for all $n \in \N_0$.\\[-6pt]

\noindent \underline{Dichotomy:}\, There are sequences $\{x_n\}_{n\in\N_0}$, $\{M_n\}_{n\in\N_0}$, $\{N_n\}_{n\in\N_0} \subset \R$ and a real number $\lambda \in (0,l)$ with the properties that $M_n, N_n \to \infty$, $M_n/N_n \to 0$,
\begin{equation}\label{eq:dichotomy}
\int_{B_{M_n}(x_n)} e_n \dx \to \lambda \quad\mbox{ and }\quad \int_{B_{N_n}(x_n)} e_n\dx \to \lambda
\end{equation}
as $n \to \infty$.
\end{theorem}

We proceed by applying Theorem~\ref{theorem:concentration} to the functions $e_n = u_n^2$, $n \in \N_0$, where $\{u_n\}_{n\in\N_0}$ is a minimising sequence for $\Efunc$ over $U_{\mu}$ with the property that $\sup_{n \in \N_0} \|u_n\|_1 < R$, so that $\ell=2\mu$.

It is a straightforward matter to exclude `vanishing'.


\begin{lemma}\label{lemma:novanishing}
No subsequence of $\{e_n\}_{n\in\N_0}$ has the `vanishing' property. 
\end{lemma}
\proof
Suppose that $\{e_n\}_{n\in\N_0}$ satisfies~\eref{eq:vanishing}, and observe that
\begin{eqnarray*}
\left| \Nfunc(u_n) \right| &\leq &\int_{\R} |N(u_n)|\dx\\ 
&\leq & c\sum_{j\in\Z} \int_{2j-1}^{2j+1}  |u_n|^{p+1} \dx\\
& \leq & c\|u\|_\infty^{p-1} \sum_{j\in\Z} \int_{2j-1}^{2j+1}  |u_n|^2 \dx\\
& \leq & c\|u\|_1^{p-1} \; \left( \sup_{x_0 \in \R} \int_{B_1(x_0)} e_n\dx \right)\\
& \leq & c  \sup_{x_0 \in \R} \int_{B_1(x_0)} e_n\dx \\
& \to 0
\end{eqnarray*}
as $n \to \infty$, which contradicts Proposition~\ref{prop:nneg}(i).\qed


\begin{lemma}\label{lemma:concentration}
Choose $s \in (0,1)$ and suppose that a subsequence of $\{e_n\}_{n\in\N_0}$ `concentrates'.
There exists a subsequence of $\{u_n(\cdot + x_n)\}_{n\in \N_0}$ which converges in $H^s(\R)$ to a minimiser of $\Efunc$ over $U_{\mu}$.
\end{lemma}
\proof Write $v_n := u_n(\cdot + x_n)$, so that $\sup_{n \in\N_0} \|v_n\|_1 < R$.
Equation ~\eref{eq:concentration} implies that for any $\varepsilon > 0$ there exists $r > 0$ such that
$$
\| v_n \|_{L^2(|x| > r)} <\varepsilon.
$$
On the other hand $\{v_n\}_{n\in\N_0}$ converges weakly in $H^1(\R)$ and strongly in $L^2(-r,r)$
to a function $v$ with $\|v\|_1 < R$; it follows that $v_n \to v$ in $L^2(\R)$ as $n \rightarrow \infty$. In view of the
interpolation inequality $\|v_n - v\|_s \leq \|v_n - v \|_0^{1-s} \|v_n - v\|_1^s$ we conclude that $v_n \to v$ in $H^s(\R)$ as $n \rightarrow \infty$, so that $\Efunc(v_n) \to \Efunc(v)$ as $n \to \infty$ (Proposition~\ref{prop:Ln2}(ii)) with $\Efunc(v)=I_\mu$ (by uniqueness of
limits).\qed


Suppose now that `dichotomy' occurs, and that $\{e_n\}_{n \in \N_0}$ satisfies~\eref{eq:dichotomy};
note in particular that the sequence $\{v_n\}_{n \in \N_0}$ with $v_n=u_n(\cdot+x_n)$ satisfies
\begin{equation}\label{eq:wn0}
\|v_n\|_{L^2(M_n < |x| < N_n)}^2 = \int_{-N_n}^{N_n} e_n \dx - \int_{-M_n}^{M_n} e_n\dx \to 0
\end{equation}
as $n \rightarrow \infty$. Let $\zeta$ be a smooth, decreasing `cut-off' function with
$$
\zeta(r) :=
\left\{\begin{array}{ll}
1, \qquad &0 \leq r \leq 1,\\[1mm]
0, \qquad &r \geq 2,
\end{array}\right.
$$
and define
\begin{eqnarray*}
v_n^{(1)}(x) &:= & v_n(x) \zeta\left( \frac{|x|}{M_n}\right),\\
v_n^{(2)}(x) &:= & v_n(x) \left( 1 - \zeta\left( \frac{2|x|}{N_n} \right) \right),
\end{eqnarray*}
so that
$$
\afl \supp ( v_n^{(1)} ) \subseteq [-2M_n, 2M_n], \qquad
\supp ( v_n^{(2)} ) \subseteq \R \setminus (-\tfrac{N_n}{2}, \tfrac{N_n}{2}),
$$
which in view of the properties of $M_n$ and $N_n$ are disjoint sets for large values of $n$.


\begin{proposition}\label{prop:wn}
The sequences $\{v_n^{(1)}\}_{n\in\N_0}$ and $\{v_n^{(2)}\}_{n\in\N_0}$ satisfy
$\|v_n^{(1)}\|_0^2 \to \lambda$, $\|v_n^{(2)}\|_0^2 \to 2\mu -\lambda$ and
\begin{equation}\label{eq:wnj0}
\|v_n^{(j)}\|_{L^2(M_n < |x| < N_n)}^2 \to 0, \qquad j= 1,2.
\end{equation}
as $n \to \infty$. 
\end{proposition}
\proof
The limits~\eref{eq:wnj0} are a direct consequence of~\eref{eq:wn0}
since $|v_n^{(1)}|$, $|v_n^{(2)}| \leq |v_n|$. It follows that
$$
\|v_n^{(1)}\|_{L^2(M_n < |x| < 2M_n)}^2 \leq \|v_n^{(1)}\|_{L^2(M_n < |x| < N_n)}^2 \to 0
$$
and
$$
\|v_n^{(1)}\|_{L^2(\frac{N_n}{2} < |x| < N_n)}^2 \leq \|v_n^{(1)}\|_{L^2(M_n < |x| < N_n)}^2 \to 0
$$
as $n \to \infty$. Using these results, we find that
\begin{eqnarray*}
\|v_n^{(1)}\|_0^2 & = & \|v_n^{(1)}\|_{L^2(|x| < M_n)}^2 + \|v_n^{(1)}\|_{L^2(M_n < |x| < 2M_n)}^2 \\
& = & \int_{-M_n}^{M_n} v_n \dx + \|v_n^{(1)}\|_{L^2(M_n < |x| < 2M_n)}^2 \\
& \to & \lambda
\end{eqnarray*}
and
\begin{eqnarray*}
\|v_n^{(2)}\|_0^2 & = & \|v_n^{(2)}\|_{L^2(|x| > N_n)}^2 + \|v_n^{(2)}\|_{L^2(\frac{N_n}{2} < |x| < N_n)}^2 \\
& = & \underbrace{\|v_n\|_2^2}_{\displaystyle = 2\mu} - \int_{-N_n}^{N_n} v_n \dx + \|v_n^{(2)}\|_{L^2(\frac{N_n}{2} < |x| < N_n)}^2 \\
& \to & 2\mu -\lambda
\end{eqnarray*}
as $n \to \infty$.\qed

Define
$$
u_n^{(1)} :=  \frac{\sqrt{\lambda}}{\|v_n^{(1)}\|_0} \, v_n^{(1)}, \qquad u_n^{(2)} :=  \frac{\sqrt{2\mu - \lambda}}{\|v_n^{(2)}\|_0} \, v_n^{(2)},
$$
so that
\begin{equation}\label{eq:unL2}
\| u_n^{(1)}\|_0^2 = \lambda, \qquad \| u_n^{(2)} \|_0^2 = 2\mu - \lambda,
\end{equation}
for all $n \in \N_0$. According to the next proposition we can assume without loss of generality that $\{u_n^{(1)}\} \subset U_\frac{\lambda}{2}$
and $\{u_n^{(2)}\} \subset U_{\mu-\frac{\lambda}{2}}$.


\begin{proposition}\label{prop:un}
The sequences $\{u_n^{(1)}\}_{n\in\N_0}$ and $\{u_n^{(2)}\}_{n\in\N_0}$ satisfy
\begin{itemize}
\item[(i)] $\lim\limits_{n \to \infty} \| v_n - u_n^{(1)} - u_n^{(2)} \|_0^2 = 0$,\\
\item[(ii)] $\limsup\limits_{n \to \infty} \|u_n^{(1)} + u_n^{(2)} \|_1  < R$ and $\limsup\limits_{n \to \infty} \|u_n^{(j)}\|_1 < R$, $j = 1,2$.\\
\end{itemize}  
\end{proposition}
\proof (i) Clearly
\begin{equation}
\| v_n - v_n^{(1)} - v_n^{(2)}\|_0^2 = \| v_n - v_n^{(1)} - v_n^{(2)}\|_{L^2(M_n < |x| < N_n)}^2 \to 0
\label{eq:wn-convergence}
\end{equation}
as $n \to \infty$ in view of the triangle inequality and the limits~\eref{eq:wn0} and~\eref{eq:wnj0}. On the other hand
\begin{eqnarray}
\afl \| u_n^{(1)} + u_n^{(2)} - v_n^{(1)} - v_n^{(2)} \|_0^2
& = & \|  u_n^{(1)} - v_n^{(1)} \|_0^2 +  \| u_n^{(2)} - v_n^{(2)} \|_0^2\nonumber \\
& = & \left( \frac{\sqrt{\lambda}}{\|v_n^{(1)}\|_0} -1 \right)^2  \| v_n^{(1)}\|_0^2 + \left( \frac{\sqrt{2\mu - \lambda}}{\|v_n^{(2)}\|_0} -1 \right)^2  \| v_n^{(2)}\|_0^2 \nonumber \\
& \to & 0
\label{eq:unwn}
\end{eqnarray}
as $n \to \infty$ (Proposition~\ref{prop:wn}).\\ 

(ii) Note that $\|v_n^\prime\|_0^2 \leq R$ and
$$\left|\frac{\mathrm{d}}{\mathrm{d}x}\! \left( \zeta\!\left(\frac{|x|}{M_n}\right) + 1 - \zeta\!\left(\frac{2|x|}{N_n}\right) \right)\right|
\leq c M_n^{-1},$$
uniformly over $x \in \R$, whence 
$$
\| ( v_n^{(1)} + v_n^{(2)} )^\prime \|_0^2 \leq \|v_n^\prime \|_0^2 + \bigO(M_n^{-1}),
$$
and~\eref{eq:wn-convergence} shows that
$$
\| v_n^{(1)} + v_n^{(2)}\|_0^2 = \|v_n\|_0^2 + o(1)
$$
as $n \to \infty$. Combining these estimates, one finds that
$$
\| v_n^{(1)} + v_n^{(2)} \|_1^2 \leq \| v_n \|_1^2 + o(1),
$$
which in the light of~\eref{eq:unwn} implies that
$$
\| u_n^{(1)} + u_n^{(2)} \|_1^2 \leq \| v_n \|_1^2 + o(1)
$$
as $n \to \infty$.

The previous inequality shows that
$$\limsup_{n \to \infty} \|u_n^{(1)} + u_n^{(2)}\|_1 \leq \sup_{n \in \N_0} \|v_n\|_1 < R,$$
and the results for $\limsup_{n \to \infty}\|u_n^{(2)}\|_1$, $j=1,2$ follow from the estimates
$$
\|u_n^{(j)}\|_1 \leq \| u_n^{(1)} + u_n^{(2)} \|_1, \qquad j =1,2.\eqno{\Box}
$$

Our next result shows that $\{\Efunc(v_n)\}_{n \in \N_0}$ decomposes into two parts for large values of $n$.

\begin{proposition}\label{prop:splitting}
The sequences $\{u_n^{(1)}\}_{n\in\N_0}$ and $\{u_n^{(2)}\}_{n\in\N_0}$ satisfy
$$\lim\limits_{n \to \infty}  \big( \Efunc(v_n) - \Efunc(u_n^{(1)}) - \Efunc(u_n^{(2)}) \big) = 0.$$
\end{proposition}
\proof First note that
\begin{equation}\label{eq:Lconvergence1}
|\Efunc(v_n)-\Efunc(u_n^{(1)} + u_n^{(2)})| \leq \sup_{u \in U} \|\Efunc^\prime(u)\|_0
\|v_n-u_n^{(1)}-u_n^{(2)}\|_0 \to 0
\end{equation}
as $n \rightarrow \infty$ since $\|\Efunc^\prime(u)\|_0$ is bounded on $U$.

Furthermore
$$
\Lfunc(u_n^{(1)} + u_n^{(2)}) = \Lfunc(u_n^{(1)}) +  \Lfunc(u_n^{(2)}) - \int_\R u_n^{(2)} L u_n^{(1)}\dx, 
$$
and
\begin{eqnarray*}
\left| \int_\R u_n^{(2)} L u_n^{(1)}\dx \right| &\leq &\tilde C_1 \|u_n^{(1)}\|_0 \int_\R \frac{|u_n^{(2)}(x)|}{\dist(x,\supp(u_n^{(1)})) }\dx \\
&\leq& \tilde C_1 R \int_{|x| > \frac{N_n}{2}} \frac{|u_n^{(2)}(x)|}{\dist(x,[-2M_n,2M_n])}\dx \\
&\leq& \tilde C_1 R^2 \left( 2 \int_{N_n/2}^\infty \frac{dx}{(x-2M_n)^2} \right)^\frac{1}{2} \\
&=& \tilde C_1 R^2 \left( \frac{4}{N_n \left( 1 - \frac{4M_n}{N_n}\right)} \right)^\frac{1}{2} \\
& \to & 0
\end{eqnarray*}
as $n \to \infty$, so that
$$
\lim_{n \rightarrow \infty} \big( \Lfunc(u_n^{(1)} + u_n^{(2)})-\Lfunc(u_n^{(1)}) -\Lfunc(u_n^{(2)}) \big) = 0.
$$
Combining this result with the equation
$$
\Nfunc (u_n^{(1)} + u_n^{(2)})  = \Nfunc ( u_n^{(1)} )  + \Nfunc (u_n^{(2)})
$$
(the supports of $u_n^{(1)}$ and $u_n^{(2)}$ are disjoint), one finds that
\begin{equation}\label{eq:Lconvergence2}
\lim_{n \rightarrow \infty} \big( \Efunc(u_n^{(1)} + u_n^{(2)})-\Efunc(u_n^{(1)}) -\Efunc(u_n^{(2)}) \big) = 0.
\end{equation}

The stated result follows from~\eref{eq:Lconvergence1} and~\eref{eq:Lconvergence2}.\qed


\begin{lemma}\label{lemma:nodichotomy}
No subsequence of $\{v_n\}_{n\in N}$ has the `dichotomy' property.
\end{lemma}
\proof
Recall that $\{v_n\}_{n \in \N_0}$ is a minimising sequence for $\Efunc$ over $U_{\mu}$ and that $\Efunc(u_n^{(1)}) \geq I_\frac{\lambda}{2}$, $\Efunc(u_n^{(2)}) \geq I_{\mu-\frac{\lambda}{2}}$. Using Lemma~\ref{prop:splitting} and the strict-subadditivity of $\mu \mapsto I_\mu$ on $(0,\mu^\star)$, we arrive at the contradiction
\begin{eqnarray*}
I_\mu &<& I_{\mu_1} + I_{\mu_2}\\ 
&\leq& \lim_{n \to \infty} \left(\Efunc(u_n^{(1)}) + \Efunc(u_n^{(2)})\right)\\ 
&=& \lim_{n \to \infty} \Efunc(v_n) \\
& = & I_\mu.\\[-6mm]
\end{eqnarray*}
\qed

According to Theorem~\ref{theorem:concentration} and 
Lemmata~\ref{lemma:novanishing} and~\ref{lemma:nodichotomy} a subsequence of
$\{e_n\}_{n \in \N_0}$ concentrates, so that the hypotheses of Lemma~\ref{lemma:concentration}
are satisfied. It follows that $D_\mu$ is nonempty and Theorem~\ref{theorem:main}(ii) holds.
The remaining assertions in Theorem~\ref{theorem:main}(i) are proved by applying
Proposition~\ref{prop:speed} and Corollary~\ref{cor:supercritical} to $u \in D_\mu$.

\section{Consequences of the existence theory} \label{sec:consequences}

\subsection*{An a priori result for supercritical solitary waves}

We now record an \emph{a priori} estimate for supercritical solutions $u \in U_\mu$ of~\eref{eq:steadyproblem}.
The result states that such solutions are long waves which `scale' in a fashion similar to the
\emph{Ansatz}~\eref{eq:longwaveansatz}. More precisely, we show that $\nn u \nn_{\tau,\mu}^2 \leq c_\tau \mu$
for $\tau<1$, where $\nn \cdot \nn_{\tau,\mu}$ is the weighted norm for $H^{2j_\star}(\R)$
defined by formula~\eref{eq:weightednorm}, so that $\|u^{(j)}\|_0 \leq c_\tau\mu^{\frac{1}{2}+j\tau\beta}$
for $j=1,\ldots,2j_\star$. We make the following additional assumption on the
nonlinearity $n$, which ensures that $u \in H^{2j_\star}(\R)$ with $\|u\|_{2j_\star} = O(\mu^\frac{1}{2})$
(see Lemma~\ref{lemma:regularity}).

{\it\begin{itemize}
\item[(A4)]
The nonlinearity $n$ belongs to $C^{2j_\star}(\R)$ with
$$n_\mathrm{r}^{(j)}(x)=\bigO(|x|^{p+\delta-j}), \qquad j=0,\ldots,2j_\star$$
for some $\delta>0$ as $x \to 0$.
\end{itemize}}

\begin{lemma}\label{lemma:scscaling}
Suppose that the additional regularity assumption (A4) holds. Every supercritical solution $u \in U_\mu$ of~\eref{eq:steadyproblem}
satisfies $\nn u \nn_{\tau,\mu}^2 \leq c_\tau \mu$ for all $\tau < 1$.
\end{lemma}
\proof
Write~\eref{eq:steadyproblem} as
$$(\nu-m)\hat{u}_1 = \xi\FF[ n(u)] , \qquad \hat{u}_2 = (\nu-m)^{-1}(1-\xi) \FF[ n(u) ].$$
Observing that
$$\nu - m(k)  > (\nu-m(0)) - \frac{cm^{(2j_\star)}(0)}{(2j_\star)!} k^{2j_\star} > - \frac{cm^{(2j_\star)}(0)}{(2j_\star)!} k^{2j_\star}$$
for $|k| < k_0$, we find that
$$
\int_\R |u_1^{(2j_\star)}|^2\dx
\leq c \int_\R ( \nu - m(k) )^2 |\hat u_1(k)|^2\dk \\
\leq c \|n(u)\|_0^2\\
\leq c  \| u \|_\infty^{2(p-1)}\|u\|_0^2.
$$
On the other hand
$$
\FF^{-1}[(\nu-m)^{-1} (1-\xi)\FF(\cdot)] \in B(L^2(\R),L^2(\R)),
$$
where the operator norm is bounded uniformly over $\nu>m(0)$, so that
\begin{eqnarray*}
\fl\lefteqn{\|u_2^{(2j_\star)}\|_0 } \\
\afl & \leq & c \|(n(u))^{(2j_\star)}\|_0 \\  
\afl & \leq & c \sum_{i=1}^{2j_\star} \|n^{(i)}(u) B_{2j_\star,i}(u^\prime,\ldots,u^{(2j_\star-i+1)})\|_0 \\
\afl & \leq & c \sum_{i=1}^{2j_\star} \|u\|_\infty^{p-i} \sum_{J_i}\|(u^\prime)^{j_1}  \ldots (u^{(2j_\star-i+1)})^{j_{2j_\star-i+1}}\|_0 \\
\afl & \leq & c \sum_{i=1}^{2j_\star} \|u\|_\infty^{p-i} \sum_{J_i}\big\|u^\prime\big\|_{L^{\frac{4j_\star}{1}}(\R)}^{j_1} \ldots
\big\|u^{(2j_\star-i+1)}\big\|_{L^{\frac{4j_\star}{2j_\star-i+1}}(\R)}^{j_{2j_\star-i+1}}\\
\afl & \leq & c \sum_{i=1}^{2j_\star} \|u\|_\infty^{p-i} \sum_{J_i}\Big(\|u\|_\infty^{1-\frac{1}{2j_\star}} \|u^{(2j_\star)}\|_0^{\frac{1}{2j_\star}}\Big)^{j_1}
\!\!\!\ldots \Big(\|u\|_\infty^{1-\frac{2j_\star-i+1}{2j_\star}} \|u^{(2j_\star)}\|_0^{\frac{2j_\star-i+1}{2j_\star}}\Big)^{j_{2j_\star-i+1}} \\
\afl & \leq & c  \|u\|_\infty^{p-1} \|u^{(2j_\star)}\|_0,
\end{eqnarray*}
where $B_{2j_\star,i}$ denote the Bell polynomials,
$$J_i=\{(j_1,\ldots,j_{2j_\star-i+1}):\ j_1+\ldots+j_{2j_\star-i+1}=i,\ j_1+2j_2+\ldots+ (2j_\star-i+1)j_{2j_\star-i+1}=2j_\star\}$$
and the generalised H\"{o}lder and Gagliardo-Nirenberg inequalities have been used
(see Hardy, Littlewood \& P\'{o}lya~\cite[Theorem 8.8]{HardyLittlewoodPolya} and Friedman~\cite[Theorem 9.3]{Friedman}).

It follows that
\begin{eqnarray}
\int_\R |u^{(2j_\star)}|^2\dx
& \leq & c \|u\|_{2j_\star}^2  \|u\|_\infty^{2(p-1)} \nonumber \\
& \leq & c \mu \|u\|_\infty^{2(p-1)} \label{eq:tauis1later}\\
& \leq & c (\mu^{1+(p-1)\tau\beta}\nn u\nn_{\tau,\mu}^{2(p-1)}) \nonumber \\
& \leq & c\left( \mu^{1+(p-1)(\tau\beta+1)}\left(\frac{\nn u\nn_{\tau,\mu}}{\mu^\frac{1}{2}}\right)^{\!\!2(p-1)}\right), \nonumber
\end{eqnarray}
and multiplying this estimate by $\mu^{-4j_\star\tau\beta}$ and adding $\int_\R u^2 \dx = 2\mu$ yields
$$
\nn u \nn_{\tau,\mu}^2 \leq c  \mu \left(1+ \mu^{(1-\tau)(p-1)}\left(\frac{\nn u\nn_{\tau,\mu}^2}{\mu}\right)^{\!\!p-1}\right).
$$
The stated estimate is obtained from this inequality using the argument given at the end of the proof of Theorem~\ref{theorem:scaling gen}.\qed

\subsection*{Convergence to long waves}

In this section we work under the additional regularity condition (A4) and
examine the relationship between $D_\mu$ and $D_\mathrm{lw}$, beginning with that
between the quantities
$$I_\mu :=\inf \left\{ \Efunc(u) \colon u \in U_{\mu} \right\}$$
and
$$
I_\mathrm{lw} := \inf \{\Efunc_\mathrm{lw}(w): w \in W_1\}.$$

\begin{lemma} \label{lemma:Iconv} \hspace{1cm}
\begin{itemize}
\item[(i)]
The quantity $I_\mu$ satisfies
$$I_\mu = -m(0)\mu+\Efunc_\mathrm{lw}(u) + o(\mu^{1+(p-1)\alpha})$$
uniformly over $u \in D_\mu$.
\item[(ii)]
The quantities $I_\mu$ and $I_\mathrm{lw}$ satisfy
$$
I_\mu = -m(0)\mu + \mu^{1+(p-1)\alpha}I_\mathrm{lw} +  o(\mu^{1+(p-1)\alpha}).
$$
\end{itemize}
\end{lemma}
\proof (i) Using the identity
$$\Efunc(u) = - m(0)\mu + \Efunc_\mathrm{lw} (u) - \frac{1}{2}\int_\R r(k)|\hat{u}|^2 \dk - \int_\R N_\mathrm{r}(u) \dx$$
for $u \in U_\mu \cap H^{j_\star}(\R)$, we find that
$$
I_\mu =  \Efunc(u)
= - m(0)\mu + \Efunc_\mathrm{lw} (u) - \frac{1}{2}\int_\R r(k)|\hat{u}|^2 \dk - \int_\R N_\mathrm{r}(u) \dx.
$$
for each $u \in D_\mu$, where
\begin{eqnarray*}
\afl \qquad\quad\left| \frac{1}{2}\int_\R r(k)|\hat{u}|^2 \dk + \int_\R N_\mathrm{r}(u) \dx \right| \\
\leq c\left( \int_\R k^{2j_\star+2} |\hat{u}|^2 \dk + \|u\|_1^2 \|u\|_\infty^{p+\delta-1}\right) \\
\leq c_\tau( \mu^{2(j_\star+1)\tau\beta} \nn u \nn_{\tau,\mu}^2 + \mu^{1+\frac{1}{2}(p+\delta-1)\tau \beta} \nn u \nn_{\tau,\mu}^{p+\delta-1}) \\
\leq c_\tau(\mu^{1+2(j_\star+1)\tau\beta} + \mu^{1+\frac{1}{2}(p+\delta-1)(\tau\beta+1)}) \\
= o(\mu^{1+(p-1)\alpha})
\end{eqnarray*}
uniformly over $u \in D_\mu$.

(ii) Choosing $u \in D_\mu$ and applying (i), one finds that
\begin{eqnarray*}
I_\mu & =  & -m(0)\mu+\Efunc_\mathrm{lw}(u) + o(\mu^{1+(p-1)\alpha}) \\
& = & -m(0)\mu+\mu^{1+(p-1)\alpha}\Efunc_\mathrm{lw}(S_\mathrm{lw}^{-1}u) + o(\mu^{1+(p-1)\alpha}) \\
& \geq & - m(0)\mu +\mu^{1+(p-1)\alpha} I_\mathrm{lw} + o(\mu^{1+(p-1)\alpha}).
\end{eqnarray*}
On the other hand, choosing $w \in D_\mathrm{lw}$ and applying Lemma~\ref{lemma:umu}, one finds that
\begin{eqnarray*}
I_\mu & \leq & \Efunc(S_\mathrm{lw}w) \\
& = & - m(0)\mu +\mu^{1+(p-1)\alpha} \Efunc_\mathrm{lw}(w) + o(\mu^{1+(p-1)\alpha}) \\
& = & - m(0)\mu +\mu^{1+(p-1)\alpha} I_\mathrm{lw} + o(\mu^{1+(p-1)\alpha}).
\end{eqnarray*}
\qed

Our main result shows how a scaling of $D_\mu$ converges to $D_\mathrm{lw}$ as $\mu \searrow 0$.

\begin{theorem}\label{theorem:solnconvergence}
The sets $D_\mu$ and $D_\mathrm{lw}$ satisfy
$$
\sup_{u \in D_\mu} \mathrm{dist}_{H^{j_\star}(\R)} \big(S_\mathrm{lw}^{-1}u , D_\mathrm{lw} \big)\rightarrow 0
$$
as $\mu \searrow 0$. 
\end{theorem}
\proof
Assume that the result is false. There exist $\varepsilon > 0$ and sequences $\{\mu_n\}_{n \in \N_0}
\subset (0,\mu_\star)$, $\{u_n\}_{n \in \N_0} \subset H^{2j_\star}(\R)$ with $u_n \in D_{\mu_n}$
such that $\lim_{n \to \infty} \mu_n = 0$ and
\begin{equation}\label{eq:doesnotconverge}
\inf_{w \in D_\mathrm{lw}} \|w_n - w \|_{j_\star} \geq \varepsilon,
\end{equation}
where $w_n(x) := \mu_n^{-\alpha} u_n(\mu_n^{-\beta}x)$. Using Lemma~\ref{lemma:Iconv}(i), one finds that
\begin{eqnarray*}
I_{\mu_n} & = & -m(0)\mu_n + \Efunc_\mathrm{lw} (u_n)+o(\mu_n^{1+(p-1)\alpha}) \\
& = & -m(0)\mu_n + \mu_n^{1+(p-1)\alpha} \Efunc_\mathrm{lw} (w_n) + o(\mu_n^{1+(p-1)\alpha})
\end{eqnarray*}
as $n \rightarrow \infty$, and because
$$I_{\mu_n} = -m(0)\mu_n + \mu_n^{1+(p-1)\alpha}I_\mathrm{lw} +  o(\mu_n^{1+(p-1)\alpha})$$
(Lemma~\ref{lemma:Iconv}(ii)), it follows that
$$\Efunc_\mathrm{lw} (w_n) = I_\mathrm{lw} + o(1)$$
as $n \to \infty$, so that $\{w_n\}_{n \in \N_0}$ is a minimising sequence for $\Efunc_\mathrm{lw}$ over
$\{w \in H^{j_\star}(\R): \Qfunc(w)=1\}$.
According to Theorem~\ref{theorem:AlbertZeng} there exists a sequence $\{x_n\}_{n \in \N_0}$ of real numbers with the property that a subsequence of $\{w_n(\cdot + x_n)\}_{n \in \N_0}$
converges in $H^{j_\star}(\R)$ to an element of
$D_\mathrm{lw}$. This fact contradicts~\eref{eq:doesnotconverge}.\qed

\begin{remark}
The previous theorem implies that $\{S_\mathrm{lw}^{-1}u\}_{u \in D_\mu}$ is a bounded set in $H^{j_\star}(\R)$.
For all $u \in D_\mu$ we therefore find that
\begin{eqnarray*}
\|u\|_\infty^2 & \leq & \frac{1}{2\pi}\|\hat u\|_{L^1(\R)}^2 \\
& \leq & \frac{1}{2\pi}\left(\int_\R\frac{1}{1+\mu^{-2j_\star\beta}k^{2j_\star}}\dk\right)\!\!
\left(\int_\R(1+\mu^{-2j_\star\beta}k^{2j_\star})|\hat{u}|^2\dk\right)\\
& = & \frac{1}{2\pi}\mu^{2\alpha}\left(\int_\R \frac{1}{1+k^{2j_\star}} \dk\right)\!\!
\left(\int_\R (1+k^{2j_\star})|\FF[S_\mathrm{lw}^{-1}u]|^2 \dk\right) \\
& \leq & c\mu^{2\alpha},
\end{eqnarray*}
and inequality~\eref{eq:tauis1later} implies that
$$\mu^{-4j_*\beta} \int_\R |u^{(2j_\star)}|^2 \dx \leq c \mu^{1+2\alpha(p-1)-4j_*\beta} = c \mu,$$
whence $\nn u \nn_{1,\mu}^2 \leq c \mu$. For $u \in D_\mu$ Lemma~\ref{lemma:scscaling}
therefore also holds with $\tau=1$ (the result predicted by the long-wave Ansatz~\eref{eq:longwaveansatz}),
and in particular $\{S_\mathrm{lw}^{-1}u\}_{u \in D_\mu}$ lies in $W$ for sufficiently large values of $S$.

\end{remark}

Finally, we relate the wave speeds $\nu(u)$ and $\nu_\mathrm{lw}(w)$
associated with respectively $u \in D_\mu$ and $w \in D_\mathrm{lw}$.

\begin{lemma} \label{lemma:wavespeedconvergence}
There exists a family $\{w_u\}_{u \in D_\mu}$ of functions
in $D_\mathrm{lw}$ such that
$$\nu(u) = m(0) + \mu^{(p-1)\alpha}\nu_\mathrm{lw}(w_u) + o(\mu^{(p-1)\alpha})$$
uniformly over $u \in D_\mu$.
\end{lemma}
\proof 
Using the identity
$$
\langle \Efunc^\prime(u),u \rangle_0 = -2m(0)\Qfunc(u)
+ \langle \Efunc_\mathrm{lw}^\prime(u),u \rangle_0 - \int_\R r(k)|\hat{u}|^2\dk
- \int_\R u\, n_\mathrm{r}(u)\dx
$$
for $u \in U$, we find that
\begin{eqnarray}
\afl\langle \Efunc^\prime(u),u \rangle_0 & = &  -2m(0)\mu
+ \mu^{1+(p-1)\alpha}\langle \Efunc_\mathrm{lw}^\prime(S_\mathrm{lw}^{-1}u),S_\mathrm{lw}^{-1}u\rangle_0 \nonumber \\
\afl & & \quad\mbox{}- \mu^{2\alpha-\beta}\int_\R r(\mu^\beta k)|\FF[S_\mathrm{lw}^{-1}u]|^2\dk
- \mu^{\alpha-\beta} \int_\R S_\mathrm{lw}^{-1}u\, n_\mathrm{r}(\mu^\alpha S_\mathrm{lw}^{-1}u)\dx \nonumber \\
\afl & = & -2m(0)\mu + \mu^{1+(p-1)\alpha}\langle \Efunc_\mathrm{lw}^\prime(S_\mathrm{lw}^{-1}u),S_\mathrm{lw}^{-1}u\rangle_0
 + o(\mu^{1+(p-1)\alpha}) \label{eq:speedconv1}
\end{eqnarray}
uniformly over $u \in D_\mu$, where the second line follows from the observation that
\begin{eqnarray*}
\afl\left|
\mu^{2\alpha-\beta}\int_\R r(\mu^\beta k)|\hat{w}|^2\dk
+\mu^{\alpha-\beta} \int_\R w\, n_\mathrm{r}(\mu^\alpha w)\dx
\right| \\
\afl\qquad\leq c\left(\mu^{2\alpha+ (2j_\star+1)\beta} \int_\R  k^{2 j_\star +2} |\hat w|^2 \dk + \mu^{(p+\delta+1)\alpha-\beta} \int_\R |w|^{p+\delta+1} \dx \right) \\
\afl\qquad = o(\mu^{1+(p-1)\alpha})
\end{eqnarray*}
uniformly over $w \in W$.

Theorem \ref{theorem:solnconvergence} asserts in particular the existence of $w_u \in D_\mathrm{lw}$
such that
$$\|S_\mathrm{lw}^{-1}u - w_u\|_{j_\star} = o(1)$$
and therefore
\begin{equation}\label{eq:speedconv2}
\afl\langle \Efunc_\mathrm{lw}^\prime(S_\mathrm{lw}^{-1}u),S_\mathrm{lw}^{-1}u\rangle_0 -
\langle \Efunc_\mathrm{lw}^\prime(w_u),w_u \rangle_0
= \sup_{ w \in W} \|\Gfunc^\prime (w)\|_0 \|S_\mathrm{lw}^{-1}u - w_u\|_0=o(1)
\end{equation}
uniformly over $u \in D_\mu$, where
$$\Gfunc(w) = \langle \Efunc_\mathrm{lw}^\prime(w),w\rangle_0
= -\int_\R\left\{\frac{m^{(2j_\star)}(0)}{(2j_\star)!} (w^{(j_\star)})^2 +(p+1) N_{p+1}(w)\right\}\dx.$$

Furthermore, it follows from the equations
$$\Efunc^\prime(u) + \nu(u) \Qfunc^\prime(u)=0, \qquad
\Efunc_\mathrm{lw}^\prime(w_u) + \nu_\mathrm{lw}(w_u)\Qfunc_\mathrm{lw}^\prime(w_u)=0$$
that
\begin{equation}\label{eq:speedconv3}
2\nu(u) \mu = - \langle \Efunc^\prime(u), u\rangle_0, \qquad
2\nu_\mathrm{lw}(w_u) = - \langle \Efunc_\mathrm{lw}^\prime(w_u), w_u\rangle_0
\end{equation}
for each $u \in U_\mu$. Combining~\eref{eq:speedconv1}--\eref{eq:speedconv3}, one finds that
$$\nu(u) = m(0) + \mu^{(p-1)\alpha}\nu_\mathrm{lw}(w_u) + o(\mu^{(p-1)\alpha})$$
uniformly over $u \in U_\mu$.\qed

\begin{remark}
For the Whitham equation Theorem~\ref{theorem:solnconvergence} and
Lemma~\ref{lemma:wavespeedconvergence} yield the convergence results
\[
\sup_{ u \in D_\mu} \inf_{y \in \R} \|\mu^{-\frac{2}{3}}  u( \mu^{-\frac{1}{3}}(\cdot+y)) -  w_\mathrm{Kdv} \|_1 \to 0
\]
and
\[
\sup_{ u \in D_\mu} \left| \nu(u) - 1 - \mu^\frac{2}{3} \left(\tfrac{2}{3}\right)^{\!\frac{1}{3}} \right| = o(\mu^\frac{2}{3})
\]
as $\mu \searrow 0$, which show how Whitham solitary waves are approximated by a scaling of the classical Korteweg-deVries solitary wave.
\end{remark}

\subsection*{Stability}

In this section we explain how Theorem~\ref{theorem:main}(ii) implies that the set of solitary-wave solutions
to~\eref{eq:DEproblem} defined by $D_\mu$
enjoys a certain type of stability, working with the following local well-posedness assumption.
(Although consideration of the initial-value problem is outside the scope of this paper we note that a local well-posedness
result in $H^s(\R)$ for $s>\frac{3}{2}$ may be obtained using Kato's method \cite{Kato75}; see also
Abdelouhab, Bona, Felland \& Saut \cite{AbdelouhabBonaFellandSaut89}.)\\
\\
{\bf Well-posedness assumption} {\it There exists a subset $M \subset U$ with the following properties.\\[-8pt]
\begin{itemize}
\item[(i)] The closure of $M \setminus D_\mu$ in $H^1(\R)$ has a non-empty intersection with $D_\mu$.\\[-8pt]
\item[(ii)] For each initial datum $u_0 \in M$  there exists a positive time $T$ and a function $u \in C([0,T],U)$ such that $u(0) = u_0$,
$$
\Efunc(u(t)) = \Efunc(u_0), \qquad \Qfunc(u(t)) = \Qfunc(u_0)
$$
for all $t \in [0,T]$ and
$$\sup_{t \in [0,T]} \|u(t)\|_1 < R.$$
\end{itemize}}


\begin{theorem}[Conditional energetic stability]\label{theorem:stability}
Choose $s \in [0,1)$. For each $\varepsilon > 0$ there exists $\delta > 0$ such that
$$
\mathrm{dist}_{H^s(\R)}(u(t),D_\mu) < \varepsilon,
$$
for all $t \in [0,T]$ whenever
$$
u \in M, \qquad \mathrm{dist}_{H^s(\R)}(u_0,D_\mu) < \delta.
$$

\end{theorem}
\proof
Assume that the result is false. There exist $\varepsilon>0$ and sequences  $\{ u_{0,n}\}_{n\in\N_0} \subset M$, $\{T_n\}_{n\in\N_0} \subset (0,\infty)$, $\{t_n\}_{n \in \N_0} \subset[0,T_n]$ and  $\{u_n\}_{n \in \N_0} \subset C([0,T_n],U)$ such that $u_n(0) = u_{0,n}$,
$$
\Efunc(u_n(t)) = \Efunc(u_{0,n}),\
\Qfunc(u_n(t)) = \Qfunc(u_{0,n}),\qquad t \in [0,T_n]
$$
and
\begin{equation}\label{eq:epsilonbound}
\mathrm{dist}_{H^s(\R)}(u_n(t_n), D_\mu) \geq \varepsilon, \qquad
\mathrm{dist}_{H^s(\R)}(u_{0,n}, D_\mu) < \frac{1}{n}.
\end{equation}
According to the last inequality there is a sequence $\{\bar u_n\}_{n\in\N_0} \subset D_\mu$ such that
\begin{equation}\label{eq:u0n-convergence}
\lim_{n \to \infty} \| u_{0,n} - \bar u_n \|_s = 0.
\end{equation}

The sequence $\{\bar{u}_n\}_{n \in \N_0}$ is clearly a minimising sequence for $\Efunc$ over $U_\mu$ with
$\sup_{n \in \N_0} \|\bar{u}_n\|_1 < R$. It follows from Theorem~\ref{theorem:main}(ii) that there is a sequence
$\{x_n\}_{n \in N_0} \subset \R$ with the property that (a subsequence of) $\{\bar{u}_n(\cdot+x_n)\}$ converges in $H^s(\R)$
to a function $\bar{u} \in D_\mu$. Equation~\eref{eq:u0n-convergence} shows that the same is true of
$\{u_{0,n}(\cdot+x_n)\}$, and using Proposition~\ref{prop:Ln2} we find that
$$\Efunc(u_{0,n}) \to \Efunc(\bar u),\qquad
\mu_n := \Qfunc(u_{0,n}) \to \Qfunc(\bar u) = \mu
$$
as $n \to \infty$. Defining $v_n := (\mu/\mu_n)^\frac{1}{2} u_n(t_n)$, observe that
$$\Qfunc(v_n)=\frac{\mu}{\mu_n}\Qfunc(u_n(t_n)) = \frac{\mu}{\mu_n}\Qfunc(u_{0,n})=\mu$$
and
\begin{eqnarray*}
\Efunc(v_n) - \overbrace{\Efunc(u_{0,n})}^{\displaystyle \to I_\mu} & = &  \Efunc(v_n) - \Efunc(u_n(t_n)) \\
& \leq & \sup_{u \in U} \| \Efunc^\prime(u) \|_0 \|v_n - u_n(t_n)\|_0 \\
& = & \sqrt{2}\sup_{u \in U} \|\Efunc^\prime(u) \|_0 |\mu-\mu_n|^\frac{1}{2} \\
& \to & 0
\end{eqnarray*}
as $n \to \infty$, so that $\{v_n\}_{n\in\N_0}$ is also a minimising sequence for $\Efunc$ over $U_\mu$ with
$\sup_{n \in \N_0} \|v_n\|_1 < R$. Theorem~\ref{theorem:main}~(ii) implies that (a subsequence of)
$\{v_n\}_{n \in \N_0}$ satisfies
$\mathrm{dist}_{H^s(\R)}(v_n, D_\mu) \to 0$ as $n \to \infty$, and since
$$\| v_n - u_n(t_n)\|_s^2 = \left(\frac{\mu}{\mu_n}-1\right) \|u_n(t_n)\|_s^2 \leq R^2\left(\frac{\mu}{\mu_n}-1\right) \to 0$$
as $n \to \infty$, we conclude that
$\mathrm{dist}_{H^s(\R)}(u_n(t_n), D_\mu) \to 0$ as $n \to \infty$. This fact contradicts~\eref{eq:epsilonbound}.\qed\pagebreak

\section*{Appendix}

Here we present a short argument demonstrating that $\Efunc_\mathrm{lw}$ is bounded below over $W_1$.

Using the Gagliardo-Nirenberg and Young inequalities, we find that
\begin{eqnarray*}
\left| \int_{\R} N_{p+1}(w)\dx \right| & \leq & c \|w\|_{L^{p+1}(\R)}^{p+1} \\
& \leq & c\|w\|_0^{(1-\theta)(p+1)} \|w\|_{j_\star}^{\theta(p+1)} \\
& \leq & c \|w\|_{j_\star}^{\frac{p-1}{2j_\star}} \\
& \leq & c_\varepsilon + c\varepsilon \|w\|_{j_\star}^2,
\end{eqnarray*}
where $\varepsilon$ is a small positive number and
$$\theta=\frac{p-1}{2j_\star(p+1)}$$
(note that $(p-1)/2j_\star < 2$ by assumption (A3)). It follows that
\begin{eqnarray*}
\Efunc_\mathrm{lw}(w) & = & \Efunc_\mathrm{lw}(w) + \Qfunc(w) - \Qfunc(w) \\
& \geq & c\|w\|_{j_\star}^2 - \left| \int_{\R} N_{p+1}(w)\dx \right| - \Qfunc(w) \\
& \geq & c \|w\|_{j_\star}^2 - c_\varepsilon \\
& \geq & -c_\varepsilon
\end{eqnarray*}
for sufficiently small values of $\varepsilon$.

\mbox{}\\[5pt]
{\bf Acknowledgement.} E. W. was supported by an Alexander von Humboldt Research Fellowship.

\bibliographystyle{nonlinearity}
\bibliography{mdg}

\section*{References}
\begin{thebibliography}{10}

\bibitem{AbdelouhabBonaFellandSaut89}
Abdelouhab L, Bona J L, Felland M and Saut J C 1989
Nonlocal models for nonlinear, dispersive waves
 {\em Physica D\/} {\bf 40} 360--392

\bibitem{Albert99}
Albert J P 1999
Concentration compactness and the stability of solitary-wave solutions to
  nonlocal equations
 {\em Contemp.\ Math.\/} {\bf 221} 1--29

\bibitem{BourdaudSickel11}
Bourdaud G and Sickel W 2011
Composition operators on function spaces with fractional order of smoothness
 {\em RIMS Kokyuroku Bessatsu\/} {\bf B26} 93--132

\bibitem{BrezisMironescu01}
Brezis H and Mironescu P 2001
{G}agliardo-{N}irenberg, composition and products in fractional {S}obolev
  spaces
 {\em J. Evol.\ Equ.\/} {\bf 1} 387--404

\bibitem{Buffoni04a}
Buffoni B 2004
Existence and conditional energetic stability of capillary-gravity solitary
  water waves by minimisation
 {\em Arch.\ Rat.\ Mech.\ Anal.\/} {\bf 173} 25--68

\bibitem{ConstantinEscher98}
Constantin A and Escher J 1998
Wave breaking for nonlinear nonlocal shallow water equations
 {\em Acta.\ Math.\/} {\bf 181} 229--243

\bibitem{Craig85}
Craig W 1985
An existence theory for water waves and the {B}oussinesq and
  {K}orteweg-de{V}ries scaling limits
 {\em Commun.\ Part.\ Diff.\ Eqns.\/} {\bf 10} 787--1003

\bibitem{EhrnstroemKalisch09}
Ehrnstr\"{o}m M and Kalisch H 2009
Traveling waves for the {W}hitham equation
 {\em Diff.\ Int.\ Eqns.\/} {\bf 22} 1193--1210

\bibitem{Friedman}
Friedman A 1969
 {\em Partial Differential Equations\/}.
New York: Holt, Rinehart and Winston

\bibitem{Gabov78}
Gabov S 1978
On {W}hitham's equation
 {\em Sov.\ Math.\ Dokl.\/} {\bf 19} 1225--1229

\bibitem{GrovesWahlen11}
Groves M D and Wahl\'{e}n E 2011
On the existence and conditional energetic stability of solitary
  gravity-capillary surface waves on deep water
 {\em J. Math.\ Fluid Mech.\/} {\bf 13} 593--627

\bibitem{HardyLittlewoodPolya}
Hardy G, Littlewood J E and P\'{o}lya G 1988
 {\em Inequalities,\/} paperback edn.
Cambridge: C.U.P

\bibitem{Kato75}
Kato T 1975
Quasi-linear equations of evolution, with applications to partial differential
  equations
In  {\em Lecture Notes in Mathematics\/ {\bf 448} --- Spectral Theory and
  Differental Equations, Dundee 1974\/}, pages 25--70 Berlin: Springer-Verlag

\bibitem{Lions84a}
Lions P L 1984
The concentration-compactness principle in the calculus of variations. {T}he
  locally compact case, part 1
 {\em Ann.\ Inst.\ Henri Poincar\'{e} Anal.\ Non Lin\'{e}aire\/} {\bf 1}
  109--145

\bibitem{NaumkinShishmarev}
Naumkin P I and Shishmarev I A 1994
 {\em Nonlinear nonlocal equations in the theory of waves\/}.
Translations of Mathematical Monographs {\bf 133} Providence, R.I.: American
  Mathematical Society

\bibitem{SchneiderWayne00}
Schneider G and Wayne C E 2000
The long-wave limit for the water wave problem. {I}. {T}he case of zero surface
  tension
 {\em Commun.\ Pure Appl.\ Math.\/} {\bf 53} 1475--1535

\bibitem{Struwe}
Struwe M 2000
 {\em Variational Methods,\/} 3rd edn.
Berlin: Springer-Verlag

\bibitem{Whitham67}
Whitham G B 1967
Variational methods and applications to water waves
 {\em Proc.\ Roy.\ Soc.\ Lond.\ A\/} {\bf 299} 6--25

\bibitem{Whitham}
Whitham G B 1974
 {\em Linear and Nonlinear Waves\/}.
New York: Wiley-Interscience

\bibitem{Zaitsev86}
Zaitsev A A 1986
Stationary {W}hitham waves and their dispersion relation
 {\em Dokl.\ Akad.\ Nauk SSSR\/} {\bf 286} 1364--1369

\bibitem{Zeng03}
Zeng L 2003
Existence and stability of solitary-wave solutions of equations of
  {B}enjamin-{B}ona-{M}ahony type
 {\em J. Diff.\ Eqns.\/} {\bf 188} 1--32

\end{thebibliography}


%
%

\end{document}